\documentclass[final,3p,times]{elsarticle}

\usepackage{graphicx}

\usepackage{amssymb}
\usepackage{amsmath}
\usepackage{subfig}
\usepackage{graphicx}
\usepackage{url}
\usepackage{algorithm2e}
\usepackage{algpseudocode}
\usepackage{multirow}
\usepackage{wasysym}
\usepackage{marvosym}
\usepackage{color, xcolor, colortbl}
\usepackage{array, booktabs}
\usepackage{xspace}
\usepackage{enumitem}

\usepackage{bm}
\definecolor{light}{rgb}{0.8,0.8,0.8}
\definecolor{medium}{rgb}{0.6,0.6,0.6}
\definecolor{dark}{rgb}{0.4,0.4,0.4}
\definecolor{darkmed}{rgb}{0.3,0.3,0.3}
\definecolor{darkest}{rgb}{0.2,0.2,0.2}
\definecolor{Black}{rgb}{0,0,0}
\definecolor{White}{rgb}{1,1,1}
\definecolor{lightpurple}{rgb}{0.78823,0.709803,0.74509}
\definecolor{lightpurpletext}{rgb}{0.788235,0.5529411,0.658823}
\definecolor{skyblue}{rgb}{0.80392,0.866666,0.92941}
\definecolor{skybluetext}{rgb}{0.61568627,0.7647058,0.913725}
\definecolor{darkgreen}{rgb}{0.3137254,0.458823,0.18431}
\definecolor{foliagegreen}{rgb}{0.188,0.415,0.105}
\definecolor{steelbluegrey}{rgb}{0.1961,0.2353,0.2392}
\definecolor{highlightblue}{rgb}{0.4078,0.6431,0.85}
\definecolor{matlabblue}{rgb}{0,0.2705,0.85}
\definecolor{darkred}{rgb}{0.8,0.1725,0}
\definecolor{fireenginered}{rgb}{0.505,0.1411,0}
\definecolor{darkpurple}{rgb}{0.6431,0.3137,0.8509}
\definecolor{gaylordpurple}{rgb}{0.416,0.204,0.549}
\definecolor{deludedorange}{rgb}{0.7409,0.4392,0}
\definecolor{darksalmon}{rgb}{0.9137,0.411,0.706}

\newcommand{\secref}[1]{Section \ref{#1}}
\newcommand{\fref}[1]{Fig.~\ref{#1}}
\newcommand{\tref}[1]{Table~\ref{#1}}

\renewcommand{\eqref}[1]{Eq.~(\ref{#1})}
\newcommand{\eref}[1]{(\ref{#1})}

\newcommand{\mat}[1]{\textrm{\textbf{#1}}}

\newcommand{\psif}{\bm{\psi} (\bm{r}, \bm{\Omega})}
\newcommand{\psifd}{\bm{\psi} (\bm{r}, \bm{\Omega}')}

\newcommand{\xten}[1]{$\times$ 10$^{\text{#1}}$}
\newcommand{\xtenm}[1]{\times \text{10}^{\text{#1}}}

\newcolumntype{a}{>{\columncolor{light}}c}
\newcolumntype{b}{>{\columncolor{skyblue}}c}

\begin{document}

\begin{frontmatter}



\title{Goal-based angular adaptivity for Boltzmann transport in the presence of ray-effects}
\author[AMCG]{S. Dargaville}
\ead{dargaville.steven@gmail.com}
\author[AWE,AMCG]{R.P. Smedley-Stevenson}
\author[AMEC,AMCG]{P.N. Smith}
\author[AMCG]{C.C. Pain}
\address[AMCG]{Applied Modelling and Computation Group, Imperial College London, SW7 2AZ, UK}
\address[AWE]{AWE, Aldermaston, Reading, RG7 4PR, UK}
\address[AMEC]{ANSWERS Software Service, Wood PLC, Kimmeridge House, Dorset Green Technology Park, Dorchester, DT2 8ZB, UK}
\begin{abstract}
Boltzmann transport problems often involve heavy streaming, where particles propagate long distance due to the dominance of advection over particle interaction. If an insufficiently refined non-rotationally invariant angular discretisation is used, there are areas of the problem where no particles will propogate. These ``ray-effects'' are problematic for goal-based error metrics with angular adaptivty, as the metrics in the pre-asymptotic region will be zero/incorrect and angular adaptivity will not occur. In this work we use low-order filtered spherical harmonics, which is rotationally invariant and hence not subject to ray-effects, to ``bootstrap'' our error metric and enable highly refined anisotropic angular adaptivity with a Haar wavelet angular discretisation. We test this on three simple problems with pure streaming where we know \textit{a priori} where refinement should occur. We show our method is robust and produces adapted angular discretisations that match the results produced by fixed refinement with either reduced runtime or a constant additional cost with angular refinement. 
\end{abstract}
\begin{keyword}
Angular adaptivity \sep Goal based \sep Ray-effects \sep Spherical harmonics \sep Filtered \sep Wavelets \sep Boltzmann transport
\end{keyword}

\end{frontmatter}
\section{Introduction}
\label{sec:Introduction}
Ray-effects are a well-known numerical artifact common in Boltzmann transport problems (aka the garden-sprinkler effect in spectral wave modelling), where an angular discretisation without rotational invariance is used and insufficient angular resolution is applied to resolve the transport caused by advection terms. For a collocation method like discrete ordinates (S$_n$) or P$^0$ FEM in angle, this causes no energy/particles to be propogated in spatial regions bounded by the geometric direction represented by each ``ray'' (analagous ray-effects are also present for high-order FEM methods in angle). The convergence of many metrics is non-monotonic within this pre-asymptotic region, especially for pointwise metrics or quantities of interest defined over small regions of the phase space. A common example is modelling the radiation dose received in a small region at some distance from a source; the dose recorded is zero until a ray aligns with the geometry and transports particles from the source.

At the same time, goal-based error metrics are finding increasing use in Boltzmann transport problems, as adaptive algorithms can help reduce the impact of the multi-dimensional phase-space (six or seven dimensions for radiation transport, five dimensions for spectral wave). Goal-based error metrics produce an adaptive algorithm that reduces the error in a functional, rather than in a globally defined norm. This is very important for problems like that described above, where the functional would be defined as the flux/dose at the end of a duct; this flux is very small in comparison to the flux near the source and hence regular adaptivity would result in near uniform refinement. 

Computing a goal-based error metric often involves multiplying both a forward and adjoint solution (or forward/adjoint residuals) and if ray-effects are present in either the forward or adjoint solutions, then for the duct problem described the resulting error metric would be zero in those regions. This means that adaptivity will not be triggered in those regions and hence a goal-based metric is useless. This is a problem faced by many goal-based error metrics in the presence of advection, e.g., in time-dependent advection problems where a functional may be zero until a wave reaches a given point in space. Even if the solution is not zero (e.g., there is a localised source), the solution and hence the error metric will miss the contribution of the distant source due to ray-effects and refine incorrectly. 

We also see this when performing Monte-Carlo simulations with variance reduction. Deterministic methods are often used to compute forward/adjoint solutions, which are combined (in a manner very similar to computing a goal-based error metric) to form a weight window used to bias sampling in the Monte-Carlo method. As such, the presence of ray-effects can make the weight windows incorrect in important regions, resulting in poor sampling. To mitigate this, in both deterministic goal-based metrics and Monte-Carlo variance reduction, a surrogate solution which does not suffer from ray-effects is often used to form error metrics/weight windows. Commonly, a diffusion equation is used, although formally a diffusion equation only results from taking the heavy-scattering limit of the Boltzmann transport equation (BTE). Naturally this means in the streaming limit our weight windows/error metrics are approximate and this impacts their effectivity. In many problems however, this approach works well and allows the use of either goal-based spatial adaptivity in a deterministic simulation, or a weight-window to perform importance sampling across different spatial regions.

Difficulties arise however, if we wish to perform deterministic goal-based angular adaptivity or importance sampling in angle (or combined space/angle adaptivity/sampling). Clearly we cannot use a diffusion equation as a surrogate, as a diffusion solution is by definition isotropic in angle and would therefore not trigger any form of anisotropic angular refinement; it could trigger uniform refinement, but we are interested in problems which require small solid angle that cannot be practically resolved with uniform angular refinement. If we wish to avoid the problems described above and build a goal-based error metric that is robust in all parameter regimes, the only recourse is to use a surrogate based on a ``true'' solution to the BTE, that is cheap to compute and is free from ray-effects. The only angular discretisation of the BTE that is free from ray-effects is a spherical harmonics ($P_n$) discretisation of the sphere as it is rotationally invariant, although this costs $\mathcal{O}(n^2)$ in angle size to compute and results in poorly conditioned linear systems in the streaming limit due to Gibbs phenomenon. Unfortunately this is the regime in which we expect heavy ray-effects to destroy the effectivity of our error metrics. 

Recently we showed \cite{Dargaville2019a} that filtered spherical harmonics \cite{McClarren2010, Radice2013} can result in well-conditioned systems in the streaming limit that are naturally free from ray-effects, as the filtering does not destroy the rotational invariance. Using low-order FP$_n$ in a vacuum modifies the expected r$^{-2}$ drop-off from a point source, causing the drop-off happen quicker than it would otherwise, but one key benefit to using FP$_n$ is that as $n \rightarrow \infty$, the solution converges to the P$_n$ solution. Importantly we also showed that in ``duct'' problems which feature pure streaming, the ray-effect free FP$_n$ discretisation meant that goal-based functionals record a poorly converged but non-zero response even at low FP$_n$ order, i.e., the FP$_n$ solution is often in the asymptotic regime for a given problem when a non-rotationally invariant (NRI) angular discretisation would not be.

It is this property we exploit in this work, with the premise being that low-order FP$_n$ solutions can serve as a ray-effect free surrogate when computing a goal-based error metric. This metric can then be used to force a NRI angular discretisation to adapt ``correctly'' in its pre-asymptotic region, when it would otherwise have an incorrect/zero response in the functional. Once the NRI angular discretisation has refined sufficiently to have all relevant streaming paths contribute to the functional (i.e., in the asymptotic region), we do not need the FP$_n$ surrogates any longer and the error metrics computed by the NRI angular discretisation suffice. 

Previously, we showed \cite{Dargaville2019} a goal-based angular adaptivity scheme based on Haar wavelets (equivalent to a P$^0$ FEM discretisation of the sphere) achieve scalable $\mathcal{O}(n)$ solutions in some heavy streaming problems and we use this as the NRI angular discretisation in this work. This paper can be considered a combination of both \cite{Dargaville2019} and \cite{Dargaville2019a}, where we allow the $\mathcal{O}(n^2)$ FP$_n$ solutions to guide our goal-based angular adapt, until the $\mathcal{O}(n)$ Haar wavelet angular adaptivity can take over and resolve the functional to high accuracy. We believe this is the first work that shows goal-based angular adaptivity in Boltzmann transport problems that is robust to ray-effects and helps lay the ground work for robust space/angle adaptivity in these mixed/hyperbolic problems. 
\section{Background}
\label{sec:Background}
In this section we outline the key features from our previous work \cite{Dargaville2019, Dargaville2019a} and introduce the discretisations that we use to build our robust error metrics. We begin with \eref{eq:bte} which is the first order, mono-energetic steady-state Boltzmann Transport Equation (BTE)
\begin{equation}
\bm{\Omega} \cdot \nabla_{\bm{r}} \psif + \Sigma_\textrm{t} \psif - S(\psif)  = S_{\textrm{e}}(\bm{r}, \bm{\Omega}),
\label{eq:bte}
\end{equation}
where $\psif$ is the angular flux in direction $\bm{\Omega}$, at spatial position $\bm{r}$. The macroscopic total cross section is $\Sigma_\textrm{t}$ with interaction/source terms given by $S(\psif)$ and external sources, $S_\textrm{e}$. We write the $\psif$ dependent source term as the typical angular scattering operator, namely
\begin{equation}
S(\psif) = \int_{\bm{\Omega}'} \Sigma_\textrm{s} (\bm{r}, \bm{\Omega}' \rightarrow \bm{\Omega}) \psifd \textrm{d}\bm{\Omega}',
\label{eq:scatter}
\end{equation}
where $\Sigma_\textrm{s}$ is the macroscopic scatter cross-sections and particles scatter from $\bm{\Omega}'$ into $\bm{\Omega}$. We now briefly describe the spatial discretisation used in this work.
\subsection{Spatial discretisation}
\label{sec:sub-grid}
The spatial discretisation we use is a sub-grid scale FEM \cite{hughes_variational_1998, hughes_multiscale_2006, candy_subgrid_2008, buchan_inner-element_2010, Dargaville2019, Dargaville2019a}, which is stable and can be considered low-memory in comparison with standard DG methods. The solution to \eref{eq:bte} is written as $\psi = \phi + \theta$, where $\phi$ and $\theta$ are the solutions on the ``coarse'' and ``fine'' scales, respectively. We then represent the coarse scale with a continuous finite-element representation, with the fine scale using a discontinuous. If we write the finite element expansions in both spaces as 
\begin{equation}
\phi(\bm{r}, \bm{\Omega}) \approx \sum_{i=1}^{\eta_N} N_i(\bm{r}) \tilde{\phi}_i(\bm{\Omega}); \qquad \theta(\bm{r}, \bm{\Omega}) \approx \sum_{i=1}^{\eta_Q} Q_i(\bm{r}) \tilde{\theta}_i(\bm{\Omega}),
\label{eq:space}
\end{equation}
with $\eta_N$ continuous basis functions, $N_i$, and $\eta_Q$ discontinuous basis functions, $Q_i$, and $\tilde{\phi}_i$ and $\tilde{\theta}_i$ containing the associated expansion coefficients, respectively. For now, it is sufficient to consider a general finite element expansion in angle with basis functions donated $G_j(\bm{\Omega})$, with a spatially varying number of basis functions $\eta_A^i$ and $\eta_D^i$ on the coarse and fine scales respectively. We then represent our expansion coefficients $\tilde{\phi}_i$ and $\tilde{\theta}_i$ in \eref{eq:space} with the space/angle expansion coefficients $\tilde{\phi}_{i,j}$ and $\tilde{\theta}_{i,j}$ and write
\begin{equation}
\tilde{\phi}_i(\bm{\Omega}) \approx \sum_{j=1}^{\eta_A^i} G_j(\bm{\Omega}) \tilde{\phi}_{i,j}; \qquad \tilde{\theta}_i(\bm{\Omega}) \approx \sum_{j=1}^{\eta_D^i} G_j(\bm{\Omega}) \tilde{\theta}_{i,j}.
\label{eq:angle}
\end{equation}
We can then apply the FEM, integrate and apply Green's theorem to recover the discretised form of \eref{eq:bte}, namely
\begin{equation}
\begin{bmatrix}
\mat{A} & \mat{B} \\
\mat{C} & \mat{D} \\
\end{bmatrix}
\begin{bmatrix}
\bm{\Phi} \\
\bm{\Theta} \\
\end{bmatrix}
=
\begin{bmatrix}
\mat{S}_{\bm{\Phi}} \\
\mat{S}_{\bm{\Theta}} \\
\end{bmatrix},
\label{eq:SGS_full}
\end{equation}
or equivalently a Schur complement
\begin{equation}
(\mat{A} - \mat{B} \mat{D}^{-1} \mat{C}) \tilde{\bm{\Phi}} = \mat{S}_{\bm{\Phi}} - \mat{B} \mat{D}^{-1} \mat{S}_{\bm{\Theta}}.
\label{eq:SGS}
\end{equation}
We have writen the vectors containing the expansion coefficients of the coarse and fine solutions as $\tilde{\bm{\Phi}}$ and $\tilde{\bm{\Theta}}$, where the discretised source terms for both scales are $\mat{S}_{\bm{\Phi}}$ and $\mat{S}_{\bm{\Theta}}$. Please see \cite{buchan_inner-element_2010} for the explicit forms of $\mat{A}, \mat{B}, \mat{C}$ and $\mat{D}$ but we should note that $\mat{A}$ and $\mat{D}$ are standard continuous and discontinuous FEM matrices, respectively. We can reconstruct the fine solution $\bm{\Theta}$ from the coarse solution by computing 
\begin{equation}
\bm{\Theta} = \mat{D}^{-1} (\mat{S}_{\bm{\Theta}} - \mat{C} \bm{\Phi}),
\label{eq:theta}
\end{equation}
and then our discrete solution is simply the addition of both the coarse and fine solutions, namely $\bm{\Psi} = \bm{\Phi} + \bm{\Theta}$ (where the coarse solution $\bm{\Phi}$ has been projected onto the fine space). 

We also make a number of modifications to the form of \eref{eq:SGS} that depend on the specific angular discretisation used (please see \cite{buchan_inner-element_2010, Goffin2014, Dargaville_2014, Goffin2015a, Buchan2016, Adigun2018, Dargaville2019, Dargaville2019a} for more details). These include approximations to $\mat{D}$ that make it easier to invert, removing coupling from \eref{eq:SGS} when adapting with FP$_n$ to increase stability and scaling down the stabilisation provided by $\mat{D}^{-1}$ in pure scattering regions to prevent locking. With these modifications we have a stable discretisation with linear growth in memory as the number of angular basis functions is increased. One of the key advantages of our sub-grid scale formulation in comparison to a standard DG formulation is that our discretisation allows us to solve for $\tilde{\bm{\Phi}}$ and then reconstruct $\tilde{\bm{\Psi}}$. As $\tilde{\bm{\Phi}}$ is on the CG mesh it is much smaller than $\tilde{\bm{\Psi}}$, which is formed on the DG mesh, allowing us flexibility in building linear solvers. In this work we use linear basis functions for both the continuous and discontinuous spatial expansions given in \eref{eq:space}.
\subsection{Angular discretisations}
\label{sec:ang_discs}
We use two different angular discretisations heavily in this work, the non-standard Haar wavelets used in \cite{Dargaville2019} and the filtered spherical harmonics used in \cite{Dargaville2019a}. We briefly discuss both discretisations here, but for more details please see \cite{Kopp1998, McClarren2010, Radice2013, Laboure2016, Laboure2016b, Frank2016}. 
\subsubsection{Haar wavelets}
\label{sec:haars}
\begin{figure}[th]
\centering
\includegraphics[width =0.25\textwidth]{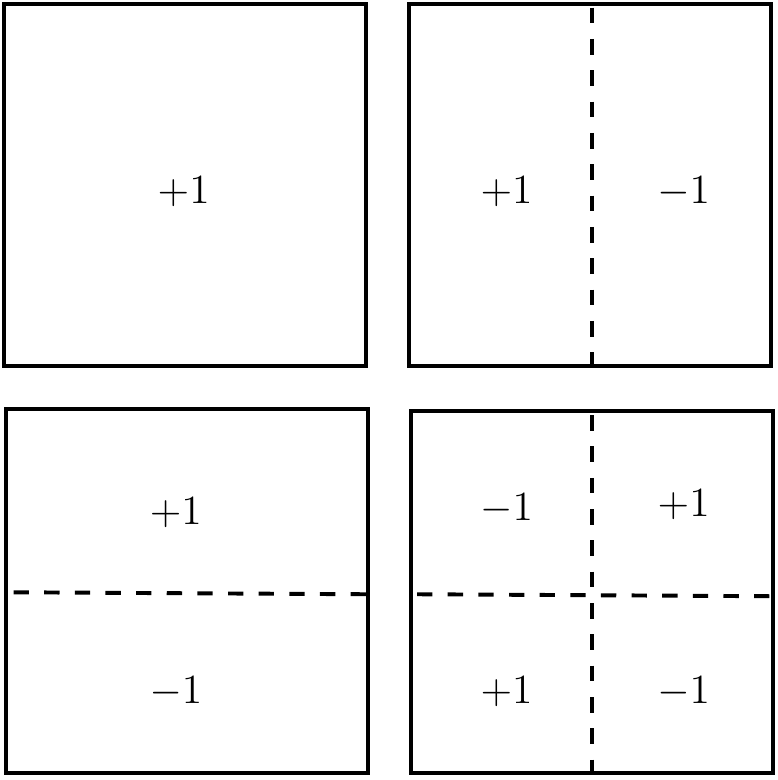}
\caption{Scaling function ($\omega_{u,n}$ - top left) and Haar wavelet basis functions ($\tau_{m,k}$) defined on an octant/quadrant, up to H$_2$. This is equivalent to discretising the octant/quadrant with a constant basis function in each of the four subpatches. }
\label{fig:wavelets}
\end{figure}

We use Haar wavelets in this work, which we build on top of a hierarchical, constant azimuthal/polar discretisation of the sphere with constant basis functions. On the two dimensional sphere, there are several equivalent Haar formulations; we use the non-standard form \cite{Kopp1998, Dargaville2019} as the wavelet basis functions have a fixed span with increasing refinement. If we consider a function $f$ on the sphere, our wavelet representation on level $j$, with $k$ wavelet functions on a given level $m$, with $n$ scaling functions on a ``base'' level $V_u$ is therefore given by
\begin{equation}
f \approx f_j = \sum_n \alpha_{u, n} \omega_{u,n} + \sum_{m=u}^{j-1} \sum_k \beta_{m,k} \tau_{m,k}, 
\label{eq:wavelet_expand}
\end{equation}
where $\alpha_{u, n}$ and $\beta_{m,k}$ are the expansion coefficients for the scaling, $\omega_{u,n}$, and wavelet functions, $\tau_{m,k}$, respectively. We use constant scaling functions on each octant/quadrant of the sphere; this gives us a base level with $n=4$ or $n=8$ and we denote each of the discretised quadrant/octant patches as i$_n$, hence our scaling functions are
\[
\omega_{u,n}(x) =     
\begin{cases}
 1, & \text{if}\ x \in i_n \\
 0, & \text{otherwise}
\end{cases}.
\]  
This is our coarsest wavelet discretisation which we refer to as H$_1$, which is roughly equivalent to an S$_2$ discretisation. Each refinement level after this then subdivides the previous into four equal area patches along constant azimuthal/polar lines taken from the halfway points in each patch. Our wavelet basis functions on each level are then patchworks of $\pm 1$ defined on the four subpatches. \fref{fig:wavelets} shows the scaling function and basis functions up to H$_2$.
Importantly our non-standard Haar wavelet discretisation is equivalent to a P$^0$ discretisation on the finest refinement level (i.e., roughly equivalent to an S$_n$ product quadrature); the wavelets trade a hierarchy in the discretisation for a hierarchy in the function expansion. This allows us to easily move between different levels of refinement with no interpolation and makes it simple to build an adaptive angular algorithm. This equivalence between the hierarchical P$^0$ space and our wavelet space also allows the use of the Mallat algorithm \cite{Mallat1989} which can map between the two spaces in $\mathcal{O}(n)$ operations. This is key to forming a scalable transport algorithm using wavelets; see \cite{Dargaville2019} for more details. 
\subsubsection{Filtered spherical harmonics}
\label{sec:fpn}
We can write a filtered version of a spherical harmonics expansion (which we refer to as FP$_n$) up to order $N$ of a function $f$, on the unit sphere as 
\begin{equation}
f(\bm{\Omega}) = \sum_{l=0}^{N} \sum_{m=-l}^{l} \left[\sigma \left( \frac{l}{N+1} \right) \right]^s f_{l,m} Y_{l,m}(\bm{\Omega}),
\label{eq:fpn}
\end{equation}
where $f_{l,m}$ are the expansion coefficients, $Y_{l,m}$ are the real, orthonormal spherical harmonics, $\sigma(\eta)$ is a filter function and $s$ is a strength parameter of the filtering. Equation \eref{eq:fpn} is rotationally invariant and if $f$ is smooth and no filter is applied ($s=0$), then \eref{eq:fpn} converges spectrally as $N$ is increased. By making $s$ nonzero, we can filter the expansion and degrade the order of convergence of our function approximation, resulting in better approximations when $f$ is not smooth.

When discretising our systems with filtered P$_n$, we use the filter $\sigma(\eta) = \sin(\eta)/\eta$ and the formulation introduced by \cite{Radice2013} which simply results in a normal P$_n$ discretisation with an extra forward-peaked scattering operator given by $-\Sigma_{\textrm{f}} \log(\sigma(l/m))$. The term $\Sigma_{\textrm{f}}$ is a free parameter (which contains $s$) which is independent of the time step/mesh spacing. We refer to this as the filter strength and this term can be made spatially dependent, allowing heavy filtering in spatial regions near discontinuities, while retaining high-order convergence in smooth regions. This results in well-conditioned linear systems, even in the presence of heavy discontinuities in space/angle and previously \cite{Dargaville2019a} we showed that the combination of FP$_n$, angular adaptivity and spatially dependent $\Sigma_{\textrm{f}}$ (computed by using the fine scale solution $\Theta$ described in \secref{sec:sub-grid} as a smoothness metric) can result in fast, accurate solutions to both smooth and non-smooth transport problems, with almost constant iteration count with angular refinement.
\section{Goal-based angular adaptivity}
\label{sec:Goal-based angular adaptivity}
We begin this section with a brief review of dual-weighted residual method we used previously \cite{Goffin2015, Dargaville2019, Dargaville2019a} to compute our goal-based error metrics. The approach we use to increase the robustness of this approach to ray-effects is detailed in \secref{sec:Adaptivity algorithm}

Our aim is to compute $\mat{e} \approx \bm{\epsilon} = \bm{\psi}_{\textrm{exact}} - \bm{\psi}$, an approximation to the exact error, $\bm{\epsilon}$, which we can use to drive our angular adaptivity. We refer to \eref{eq:bte} as the ``forward'' problem, with exact solution $\bm{\psi}_{\textrm{exact}}$ and residual $\mathcal{R}$, hence $\mathcal{R}(\bm{\psi}_{\textrm{exact}}) = 0$. We begin by considering a functional, $F$ of the solution which we are trying to minimise the error in
\[
F(\bm{\psi}) = \int_P f(\bm{\psi}) \, \textrm{d}P,
\]
where $f$ is an arbritrary function of the solution, such as the average flux in a region or current over the surface, with $P$ representing the phase-space.
\subsection{Error metric}
\label{sec:Error metric}
Following \cite{Goffin2015}, we can approximate the error in our functional as
\begin{equation}
|F(\bm{\psi}_{\textrm{exact}}) - F(\bm{\psi})| \approx \bm{\epsilon}^{\textrm{T}} \mat{R}^*
\label{eq:error}
\end{equation}
or equivalently
\begin{equation}
|F(\bm{\psi}_{\textrm{exact}}) - F(\bm{\psi})| \approx \bm{\epsilon}^{*\textrm{T}} \mat{R}
\label{eq:adjoint_err}
\end{equation}
where $\bm{\epsilon}^{\textrm{T}}$ and $\bm{\epsilon}^{*\textrm{T}}$ are the discrete forward and adjoint solution error, respectively, with $\mat{R}$ and $\mat{R}^*$ the discrete forward and adjoint residuals computed using $\bm{\psi}^*$ and $\bm{\psi}^*_{\textrm{exact}}$, which are the approximate and exact solutions of the adjoint equation with source term derived from the response function respectively.

Due to Galerkin orthogonality, we must modify the forward and adjoint residuals so they are non-zero; we do this by computing ``reduced-accuracy'' residuals, $\hat{\mat{R}}$ and $\hat{\mat{R}}^*$. We chose to reduce the accuracy with which we compute these residuals to ensure we don't introduce artifacts into our error metrics when we have highly anisotropic angular flux; many authors try to use high order interpolation to form these which can lead to the introduction of oscillations. We also must chose a target error for our adaptivity, $\tau$, and hence form our error metric by taking the pointwise maximum of both \eref{eq:error} and \eref{eq:adjoint_err} over each angular coefficient. Our forward and adjoint problems share the same adapted angular discretisation and the form of \eref{eq:gb_metric} means that features in both the forward and adjoint solutions cause angular refinement. Our approximate goal-based error metric is therefore given by
\begin{equation}
\mat{e} = \frac{\max\{|\bm{\epsilon} \odot \hat{\mat{R}}^*|, |\bm{\epsilon}^* \odot \hat{\mat{R}}|\} N_{\textrm{DOF}}}{\tau}, 
\label{eq:gb_metric}
\end{equation}
where $\odot$ denotes pointwise multiplication. Given our wavelet discretisation is hierarchical, we choose to approximate our solutions errors as $\bm{\epsilon} \approx \bm{\psi}$ and $\bm{\epsilon}^* \approx \bm{\psi}^*$; we should note the absolute values in these expressions have been removed compared with \cite{Dargaville2019}.

To compute our reduced accuracy residuals, we compute residuals on both our coarse and fine scales, $\hat{\mat{R}}_{\bm{\Phi}}$ and $\hat{\mat{R}}_{\bm{\Theta}}$ respectively, using
\begin{equation}
\begin{bmatrix}
\hat{\mat{R}}_{\bm{\Phi}} \\
\hat{\mat{R}}_{\bm{\Theta}} \\
\end{bmatrix}
=
\begin{bmatrix}
\tilde{\mat{A}} & \tilde{\mat{B}} \\
\tilde{\mat{C}} & \tilde{\mat{D}} \\
\end{bmatrix}
\begin{bmatrix}
\bm{\Phi}\\
\bm{\Theta} \\
\end{bmatrix} - 
\begin{bmatrix}
\mat{S}_{\bm{\Phi}} \\
\mat{S}_{\bm{\Theta}} \\
\end{bmatrix}.
\label{eq:disc_resid_subgrid}
\end{equation}
The modified submatrices, $\tilde{\mat{A}}, \tilde{\mat{B}}, \tilde{\mat{C}}$ and $\tilde{\mat{D}}$ are then chosen depending on the angular discretisation we are using. For the Haar wavelets, we follow \cite{Dargaville2019} and set $\mat{S}_{\bm{\Phi}}=\mat{S}_{\bm{\Theta}}=0$ and use only the diagonals of $\tilde{\mat{A}}$, $\tilde{\mat{B}}$, $\tilde{\mat{C}}$ and $\tilde{\mat{D}}$ when computing our residual.This residual can be calculated easily with a single matrix-vector product and is non-zero. The coarse and fine residuals are then combined like the discrete solution in \secref{sec:sub-grid} to form our reduced accuracy discrete residual, $\hat{\mat{R}}$, as
\begin{equation}
\hat{\mat{R}} = \hat{\mat{R}}_{\bm{\Phi}} + \hat{\mat{R}}_{\bm{\Theta}}  
\label{eq:disc_resid}
\end{equation}
\subsection{Ray-effects}
\label{sec:Ray-effects}
\begin{figure}[th]
\centering
\includegraphics[width =0.6\textwidth]{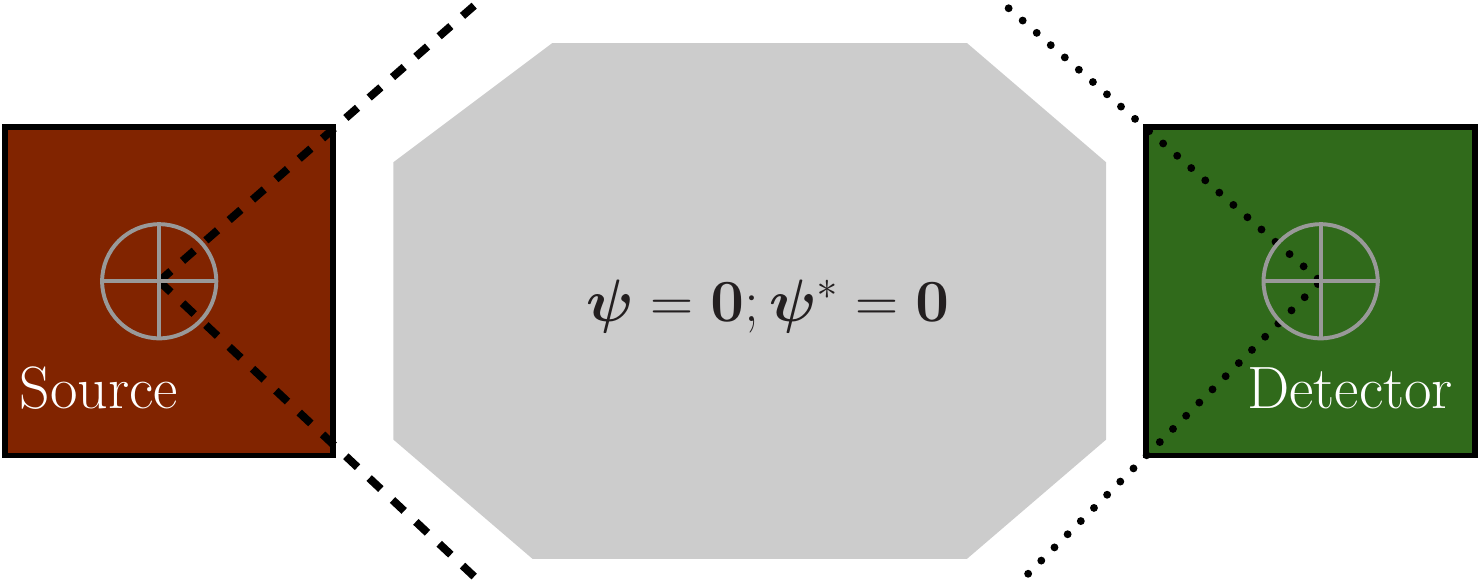}
\caption{Schematic of a source/detctor problem in a vacuum with a coarse angular discretisation. Both the solution of the forward (dashed lines) and adjoint (dotted lines) problems suffer from ray effects; the shaded grey region has zero flux in both the forward and adjoint solutions.}
\label{fig:ray_effect}
\end{figure}

As mentioned in \secref{sec:Introduction}, we must take care combining goal-based angular adaptivity with any NRI angular discretisation, as the metric described in \secref{sec:Error metric} will be incorrect/zero if ray-effects are present. If we take the example of a duct problem in a vacuum (see \fref{fig:ray_effect}), this is clear as both $\bm{\epsilon}$ and $\bm{\epsilon}^*$ (and the residuals) in \eref{eq:gb_metric} will be zero as $\bm{\psi}$ and $\bm{\psi}^*$ are zero. Our error metric, $\mat{e}$, is therefore zero.

The ``answer'' to this problem is the same required to resolve ray-effects in the solution, namely refine the angular mesh until we have sufficient resolution to partially resolve all streaming paths. This will ensure that the error metric will ``see'' all streaming paths and angular adaptivity can then refine where necessary to produce solutions with lower error. Of course this is undesirable, given that unless we know \textit{a priori} where this initial refinement must take place we are forced to refine the angular domain uniformly and the problems we are interested in solving cannot be feasibly solved using uniform angular resolution.

Previously, the only authors to investigate using goal-based angular adaptivity with NRI angular discretisations are \cite{Bennison2014, Murphy2015, Goffin2015, Goffin2015a, Hall2017, Soucasse2017, Zhang2018, Dargaville2019}, and they do not discuss any measures used to mitigate the problems discussed in this work. This is because the successful use of angular adaptivity in those works rely on four factors; solving streaming problems that do not use pure vacuums or very smooth problems without significant ray-effects, relying on scattering to get a detector response, deliberate alignment of the coarsest angular discretisation and the problem geometry, and numerical diffusion caused by insufficiently resolved spatial meshes. 

In \cite{Goffin2015, Goffin2015a} for example, the 3D duct problem modelled has a ``near'' vacuum ($\Sigma_{\textrm{t}}=0.001$) in the duct region and a scattering region surrounding the entire duct. This means that particles have a fixed propogation distance and even with the coarsest angular resolution, scattering down the entire length of the duct region ensures the detector receives a signal from the source. 

In contrast to this, we previously \cite{Dargaville2019} solved a 2D dogleg duct problem, which did not have any scattering regions and featured a pure vacuum, with the detector not aligned with streaming paths from the source when using the coarse angular discretisation (H$_1$ which is roughly equivalent to S$_2$). In that work we ``ensure[d] that any goal-based problem we run has a non-zero response even with a coarse angular discretisation''. This was achieved purely with numerical diffusion from our spatial mesh; we used only 2824 elements to resolve this problem. 
\begin{figure}[th]
\centering
\subfloat[][Coarse spatial mesh from \cite{Dargaville2019} with 2824 elements and a non-standard Haar wavelets discretisation in angle with level 1 refinement (very similar to S$_2$). The average flux in the goal region is -0.54\xten{-3}.]{\label{fig:coarse_2D_duct}\includegraphics[width =0.3\textwidth]{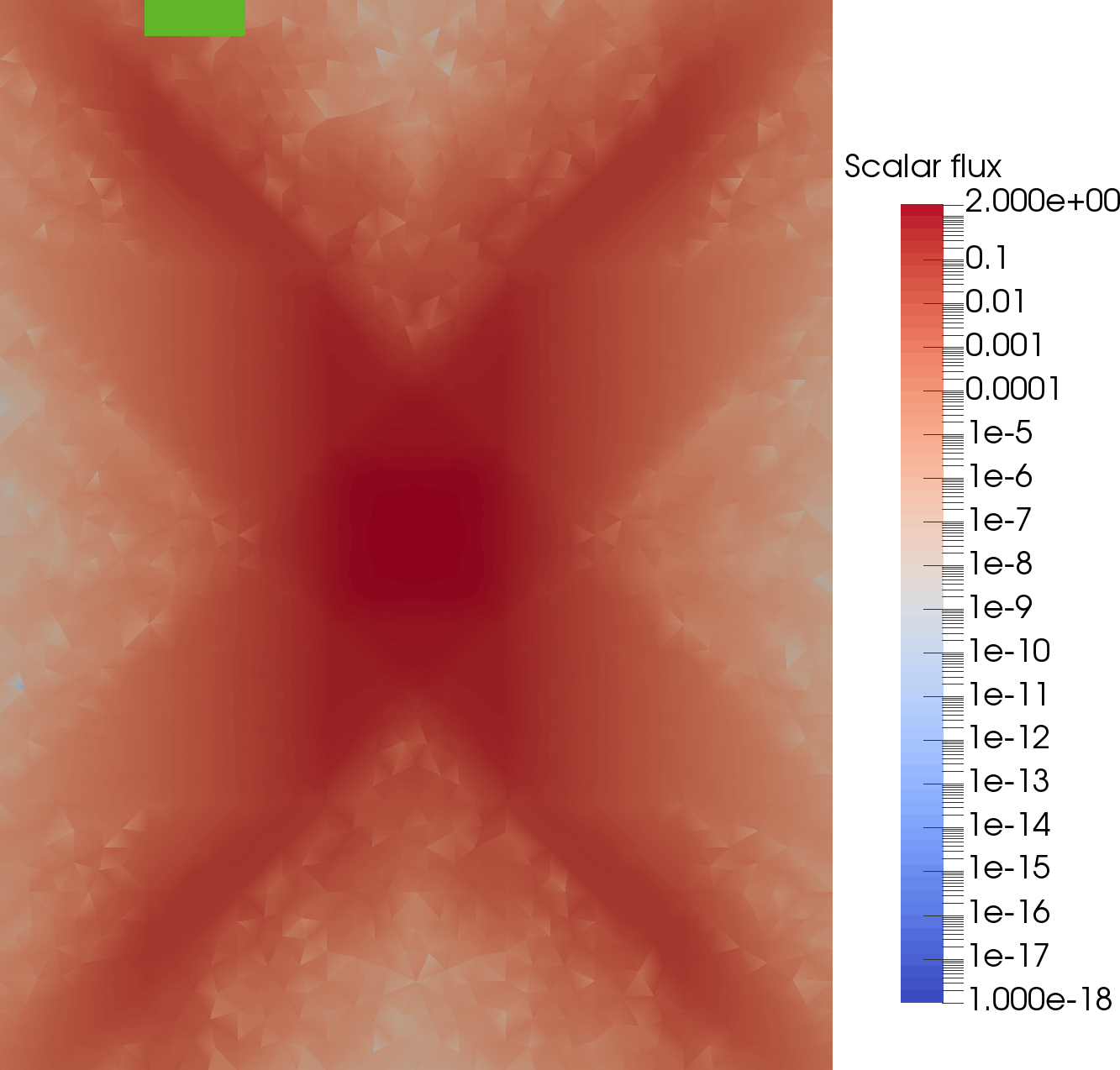}} \hspace{0.2cm}
\subfloat[][Refined spatial mesh with 264,006 elements and a non-standard Haar wavelets discretisation in angle with level 1 refinement (very similar to S$_2$). The average flux in the goal region is -0.49\xten{-10}.]{\label{fig:refine_2D_duct}\includegraphics[width =0.3\textwidth]{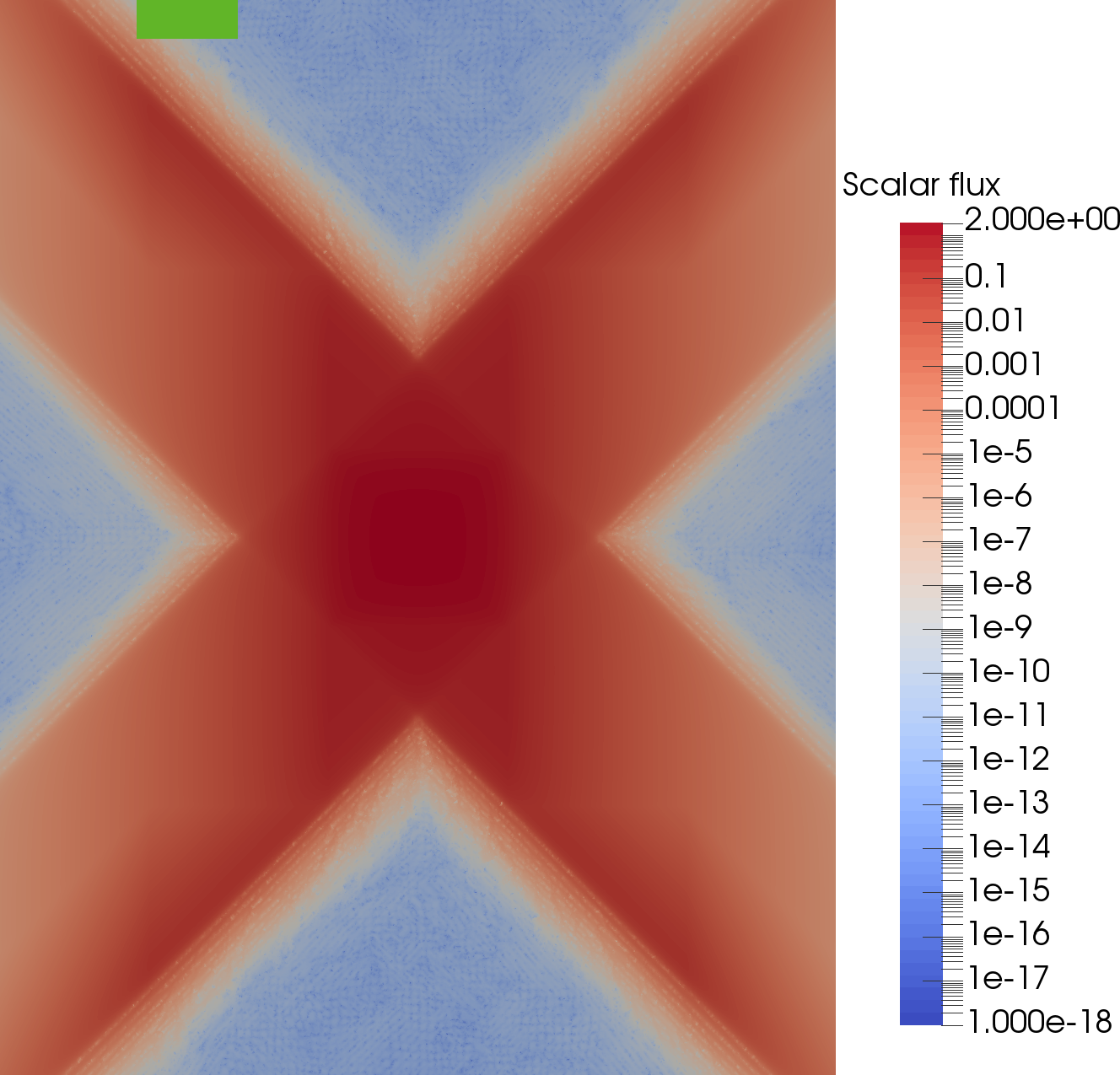}}\\
\subfloat[][Refined spatial mesh with 264,006 elements and an FP$_1$ discretisation with $\Sigma_\textrm{f}=1$. The average flux in the goal region is 0.19\xten{-1}.]{\label{fig:refine_fpn}\includegraphics[width =0.3\textwidth]{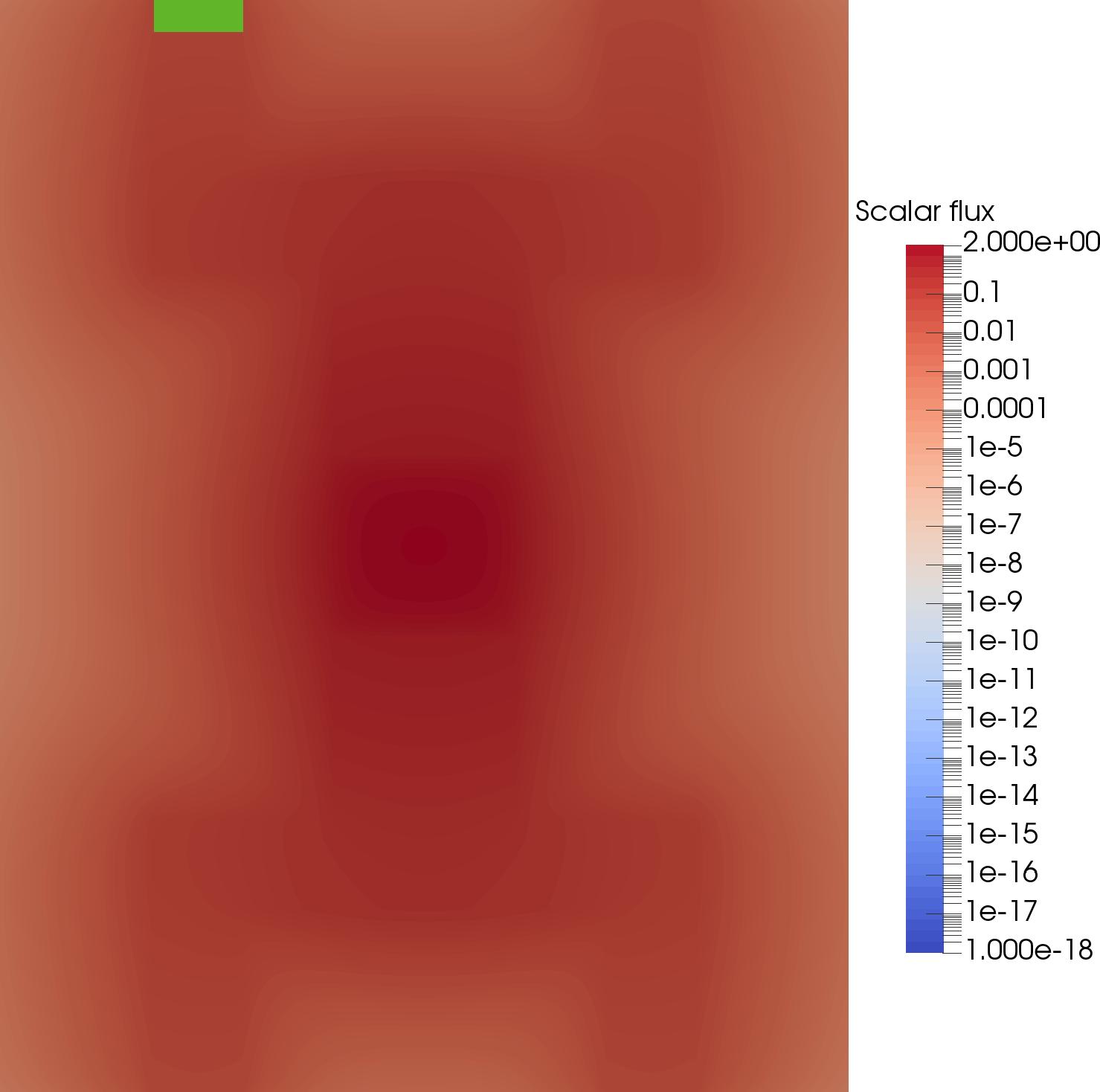}}
\caption{Absolute value of the scalar flux across the domain on a coarse mesh for the 2D dogleg problem used in \cite{Dargaville2019}. The goal-region is shaded green at the end of the duct in the top-left corner.}
\label{fig:coarse_vs_refine}
\end{figure}

This was a deliberate choice given that in that work we were trying to show the scalability of our wavelet-based angular adaptivity and could not feasibly converge both the spatial and angular meshes. If however we converge the spatial mesh in that problem (see \fref{fig:coarse_vs_refine}), we can see that with a mesh of 264k elements the detector response in that problem with a coarse angular discretisation is 10$^{-10}$ (and approaching zero with spatial mesh refinemnt, as it should) and angular refinement would not be triggered. We must not confuse the triggering of angular adaptivity in those works to mean the goal-based error metrics used are robust in the presence of ray-effects.  

Besides using a rotationally-invariant angular discretisation to guide the adaptivity as we do in this work, there are several ways we could attempt to resolve this problem, namely:
\begin{enumerate}
\item Rotate the angular quadrature/problem geometry to align with streaming paths. This is only possible in simple problems where streaming paths are known \textit{a priori}, which precludes complex multi-group problems that may include different physics (e.g. in charged particle transport with electric fields that cause turning).
\item Compute an additional forward/adjoint solution with a different NRI angular discretisation/quadrature and use that to help form the error metric. The success of this approach would rely on the ray-effects produced with different discretisations never overlapping.
\item Introduce angular diffusion/filtering to our NRI angular discretisation to smear out ray-effects (e.g., \cite{Booij1987, Tolman2002} or the recent filtered S$_n$ work of \cite{Hauck2019}). The smoothing provided will however only be effective out to a given spatial distance, as one could not add infinite angular diffusion to smooth out all ray-effects without triggering uniform angular refinement.
\item Use regular angular adaptivity (i.e., not goal-based) to produce an adapted solution that can then be used as the intial step in a goal-based adapt. For problems with heavy streaming, a very low adaptive tolerance would be required to adapt down a duct, for example, which would result in (close to) uniform refinement in angle. 
\end{enumerate}

In the next section, we outline our approach which combines both FP$_n$ and an adaptive wavelet scheme in an attempt to overcome the disadvantages faced by the methods listed above. As mentioned this is based on using an FP$_n$ solution to identify where our error metric would be underresolved in Haar space. \fref{fig:refine_fpn} for example shows the FP$_1$ solution in the 2D dogleg problem with a refined spatial mesh and we can see both the absence of ray-effects and also evidence of an non-zero detector response. 

We should note that although we can consider FP$_n$ as a P$_n$ method with additional angular diffusion (added by the ``negative'' forward peaked scatter term introduced by \cite{Radice2013} and discussed in \secref{sec:fpn}), the use of FP$_n$ in this work is not designed to tackle the lack of differentiability in the dual problem (e.g., see \cite{Beckers2019}) given our mixed hyperbolic transport problems. That is a fundamental problem which underlies the dual-weighted error metrics we use. Often, extra diffusion in the form of an artificial viscosity is added to ensure enough continuity exists for the dual problem (e.g., see \cite{Johnson1995, Pierce2004a})

Our goal in this paper is to first tackle the impact the pre-asymptotic region of any NRI angular discretisation has on goal-based error metrics. Ray-effects are numerical artifacts which introduce artificial discontinuities in our solutions when our solutions are underresolved, in contrast to the real discontinuities which exist in our solution (i.e., when our discretisation is in the asymptotic region) and destroy the continuity required by the formulation of the dual problem. The fact we use an angular discretisation with additional angular diffusion in this work to help tackle the impact of ray-effects in our error metric is convenient, as it forms a natural way to improve our error metric once our angular discretisation is in the asymptotic regime; we leave investigating this to future work. 
\subsection{Adaptivity algorithm}
\label{sec:Adaptivity algorithm}
The first step in our adaptive process involves computing a forward and adjoint solution with our Haar wavelet discretisation, with coarse resolution; in this work our first adapt step is always H$_1$. We then compute a forward and adjoint solution using a low-order FP$_n$ discretisation. There are a number of ways we could combine these four solutions to produce a robust error metrics. One option would be to simply compute a goal-based error metric using the theory described in \secref{sec:Error metric} for both the wavelets and the FP$_n$. The metric produced by the wavelets would be zero/too small in regions where ray-effects prevent contributions to the functional, whereas the FP$_n$ metric would not be. At each spatial node we could therefore project the FP$_n$ metric into Haar space and take the maximum across both metrics over each wavelet coefficient in Haar space. Naturally we would need to take care to do our thresholding in Haar space, as the thresholding we use to adapt our FP$_n$ discretisation in \cite{Dargaville2019a} only triggers refinement on each spatial node, as a FP$_n$ discretisation cannot adapt anisotropically in angle. Our thresholded error metric would then drive our anisotropic angular adapt in Haar space. We could then compute the forward and adjoint solutions with our newly adapted wavelet discretisation and repeat this process, using the same low-order FP$_n$ solution we computed in the first step.

The main problem with this approach is that it relies on the effectivity indices of both the Haars and the FP$_n$ being at best ``good' and at worst proportional to each other. Previously we showed that the effectivity index for our wavelet discretisation was pathologically bad, reaching 10$^7$ in some problems \cite{Dargaville2019}, whereas our FP$_n$ effectivity index in the same problem, while still needing improvement, performed much better \cite{Dargaville2019a}, with values between 0.4 and 22. \cite{Dargaville2019, Dargaville2019a} however showed that these poor metrics are still good enough to produce refinement where required to reduce error in our goal-based functionals. Furthermore, the route to improving these indices is clear, by computing a more sophisticated reduced accuracy residual with Haars (i.e., not using a simple diagonal matvec), ensuring our spatially-dependent forward and adjoint FP$_n$ solutions are projected to a common space, and possibly using additional diffusion to introduce the required continuity in the dual problem, as suggested in \secref{sec:Ray-effects}. 

If we were to take the maximum of the Haar and FP$_n$ error metrics as described above, it is easy to imagine cases where either of the metrics are erroneously smaller than the other, causing incorrect refinement or stagnation in the adapt process. Until we develop both a Haar and FP$_n$ error metric with ``good'' effectivity indices and are convinced that the only time our Haar effectivity index is ``wrong'' is when ray-effects are affecting the error metric, we must rely on a different approach. 

Given we know both the Haar and FP$_n$ solutions converge to the true solution of the BTE, we instead rely on the FP$_n$ solution to tell us when the Haar solution is underresolved. We use a simple heuristic to chose nodes where the FP$_n$ forward/adjoint solution is more representative of the true solution; in this work if the absolute value of either the forward/adjoint Haar solution is not within ten times of the absolute value of the forward/adjoint FP$_n$ solution on a given spatial node, we denote that node as underresolved. On each underresolved node, we project both the forward and adjoint FP$_n$ solutions into adapted Haar space. If the heuristic determines the node is sufficiently resolved, then the Haar solutions are used on that node instead. These combined FP$_n$/Haar solutions for the forward and adjoint problems are then used when computing $\bm{\epsilon}$, $\bm{\epsilon}^*$, $\hat{\mat{R}}$ and $\hat{\mat{R}}^*$. The error metric is then computed using \eref{eq:gb_metric}. 

This projection from FP$_n$ to adapted Haar space will be $\mathcal{O}(n^2)$ given the support of each FP$_n$ basis function, but computing our FP$_n$ solution is already $\mathcal{O}(n^2)$ and hence this does not change our reliance on low-order FP$_n$ solutions. We should also note that given the different approaches to computing reduced accuracy residuals in Haar space as described in \secref{sec:Error metric} and in FP$_n$ space as described in \cite{Dargaville2019a}, this is not equivalent to simply using the FP$_n$ error metric on underresolved nodes (though it is similar). 

There are two obvious drawbacks to this approach, namely
\begin{enumerate}
\item We require the FP$_n$ solution to be in the asymptotic regime.
\item We must chose a ratio of solutions that is large enough to not require fully converged solutions in the FP$_n$ discretisation and also small enough to differentiate between multiple distant sources; e.g., a scattering source also contributing some small fraction of a detector response from a distance. 
\end{enumerate}

The first of these is central to this work, the assumption that low-order FP$_n$ is in the asymptotic regime before a NRI discretisation. Thankfully this is easily quantified; our previous work \cite{Dargaville2019a} for example showed a non-zero detector response with FP$_9$ and constant filter strength $\Sigma_\textrm{f}=1$ in a pure vacuum duct with width/length ratio of 1/100. An FP$_9$ solution would therefore be sufficient to drive our adapt in all the examples shown in this work and is cheap to compute. If we require greater resolution to reach the the asymptotic FP$_n$ regime, this can become more costly. Thankfully again we benefit from our previous work \cite{Dargaville2019a} where we showed that a spatially-dependent filter strength and goal-based angular adaptivity for FP$_n$ can significantly reduce the runtime in streaming problems and help delay the onset of the $\mathcal{O}(n^2)$ behaviour. Given we are solving both a ``coarse'' forward/adjoint uniform FP$_n$ problem in this work and are stepping through an adapt process with our wavelets, it is simple to enable the FP$_n$ adaptivity and spatially-dependent filter at the same time. For simplicity we do not show this here, but we have found this significantly expands the range of problems we can feasibly tackle, with both our Haar wavelet and FP$_n$ discretisations performing goal-based angular adaptivity simultaneously. 

The second drawback mentioned above is that our algorithm is sensitive to the choice of ratio which determines ``underresolved'' nodes. This is an unfortunate consequence of being unable to rely on the error metrics in both spaces having good effectivity indices. We must note however, that as we decrease this ratio, increase the order of our FP$_n$ solution and decrease $\tau$ (the thresholding tolerance we use in our adapt), we are at least guaranteed to adapt in the presence of ray-effects, as we know the FP$_n$ solution converges to the true solution of the BTE. Indeed an easy test for the convergence of our adaptive process is to increase the FP$_n$ order used (or turn on FP$_n$ goal based adaptivity as mentioned above). This is a significant advantage when compared to goal-based error metrics which do not use some sort of surrogate; they would not adapt at all, regardless of the value of $\tau$ used. To summarise, for a generic NRI discretisation (we use our non-standard Haar wavelets) our algorithm is given by Algorithm \ref{alg:one}.

\begin{algorithm}[ht]
 Compute coarse FP$_n$ forward, $\bm{\Psi}_i^\textrm{FP$_n$}$, and adjoint, $\bm{\Psi}_i^{\textrm{FP$_n$}*}$, solutions\;
 Compute coarse NRI discretised forward, $\bm{\Psi}_i^\textrm{NRI}$, and adjoint, $\bm{\Psi}_i^{\textrm{NRI}*}$, solutions\;
 \While{the current refinement is below the defined maximum level}{
  \For{each node in the fine DG scale, $i$}{
  	\eIf{$\bm{\Psi}_i^\textrm{NRI} > 10 \bm{\Psi}_i^\textrm{FP$_n$}$ or $\bm{\Psi}_i^{\textrm{NRI}*} > 10 \bm{\Psi}_i^{\textrm{FP$_n$}*}$}{
   	$\bm{\Psi}_i^\textrm{resolved} = \bm{\Psi}_i^\textrm{FP$_n$}$ projected into NRI space\;
   	$\bm{\Psi}_i^{\textrm{resolved}*}  = \bm{\Psi}_i^{\textrm{FP$_n$}*}$  projected into NRI space\;   	
   	Similarly form $\bm{\Theta}_i^\textrm{resolved}$ and $\bm{\Theta}_i^\textrm{{resolved}*}$\;
   	Similarly form $\bm{\Phi}_i^\textrm{resolved}$ and $\bm{\Phi}_i^\textrm{{resolved}*}$ for the co-located CG node of $i$\;
   	}{
   	$\bm{\Psi}_i^\textrm{resolved} = \bm{\Psi}_i^\textrm{NRI}$\;
   	$\bm{\Psi}_i^{\textrm{resolved}*} = \bm{\Psi}_i^{\textrm{NRI}*}$\;
   	Similarly form $\bm{\Theta}_i^\textrm{resolved}$ and $\bm{\Theta}_i^\textrm{{resolved}*}$\;
   	Similarly form $\bm{\Phi}_i^\textrm{resolved}$ and $\bm{\Phi}_i^\textrm{{resolved}*}$ for the co-located CG node of $i$\;  	
  }
  }  
  Construct $\bm{\epsilon}$ and $\bm{\epsilon}^*$ using $\bm{\Psi}_i^\textrm{resolved}$ and $\bm{\Psi}_i^{\textrm{resolved}*}$\;
  Construct $\hat{\mat{R}}$ and $\hat{\mat{R}}^*$ using $\bm{\Theta}_i^\textrm{resolved}$, $\bm{\Theta}_i^{\textrm{resolved}*}$, $\bm{\Phi}_i^\textrm{resolved}$ and $\bm{\Phi}_i^{\textrm{resolved}*}$ \;
  Compute the error metric using \eref{eq:gb_metric}\;
  Refine/coarsen the NRI discretisation\;
 }
 \caption{Our adaptivity algorithm with a generic NRI angular discretisation and surrogate FP$_n$ solutions.}
 \label{alg:one}
\end{algorithm}

The remaining details of our algorithm, including how our Haar wavelets are targetted for refinement/coarsening with thresholding and our iterative method are unchanged from \cite{Dargaville2019, Dargaville2019a}. We should also note that as discussed in \cite{Dargaville2019, Dargaville2019a}, we define our adjoint angular domain as the negative of the forward problem. Given our error metric is taken from combining the forward and adjoint error measures, this allows us to use the same adapted angular domain in both the forward and adjoint problems. This means that for the results shown below, the adapted forward angular flux plots also show where adaptivity has occured in the (reflected) adjoint problem. 
\section{Results}
\label{sec:Results}
Outlined below are three example problems we use to test the robustness of our new error metric; in one problem we compare results against the error metric used in \cite{Dargaville2019}; to differentiate between them we refer to the error metric in this work as ``robust'', and the one used in \cite{Dargaville2019} and defined in \secref{sec:Error metric} as ``non-robust''. The three examples have all been chosen to be simple enough to know \textit{a priori} where angular refinement should occur to produce accurate solutions, so we can verify that our algorithm is adapting correctly. These ``fixed refinement'' simulations refine the angular domain between azimuthal and polar bounds up to some maximum refinement level, using the same angular discretisation across all spatial nodes. Solving the fixed refinement problems therefore only requires the solution of a single linear system. 

We should emphasise that we do not expect our adaptive method to produce solutions faster or with fewer DOFs in these problems when compared to fixed refinement; the goal in this work is to confirm our method is robust when tested on problems which feature heavy ray-effects. Importantly we have been careful to construct these problems so that goal-based angular refinement is not triggered correctly at coarse resolution when using the non-robust metric. As discussed in \secref{sec:Ray-effects} this means using problems which feature pure streaming, with no/small scattering regions, ensuring no alignment between the H$_1$ coarse angular discretisation (which discretises the sphere into quadrants/octants) and a detector, and using spatial meshes fine enough that numerical diffusion does not trigger a response. 

We do not solve the adjoint problem at the last adapt step as it is unecessary. All the runtimes shown for adapted simulations include all the adapt steps, computation of error metrics, etc. This includes the cost of computing the coarse FP$_n$ solutions. For example, if we take 10 adapt steps, the runtime includes the cost of solving 21 linear systems (10 forward Haar, 9 adjoint Haar, 1 forward FP$_n$ and 1 adjoint FP$_n$). Unless otherwise noted, all linear systems are solved to an absolute and relative tolerance of 1\xten{-10}. In some of the problems below we take an extra adapt step at the maximum order, to allow the adaptivity to ``settle-down'', as in \cite{Dargaville2019}.
\subsection{2D void problem}
\label{sec:2D void problem}
The first example problem is a simple 2D source/detector problem featuring pure vacuum, like that shown in \fref{fig:ray_effect}, with the source (strength 1) and detector regions of size 1x1cm separated by a distance of 10cm (so the total length of the domain is 12cm). We discretise this problem in space with an unstructured triangular mesh with 502 elements (304 CG nodes and 1506 DG nodes). We can \textit{a priori} determine that fixed angular refinement between $\mu \in [0, 1]$ and $\omega \in [1.47976, 1.661832]$ in this problem results in a discretisation that can resolve the transport of particles down the duct; we perform fixed refinement in these bounds with our non-standard Haar wavelets as a baseline to compare against in this problem. The reference solution used to compute the relative error is our non-standard Haar discretisation with fixed refinement up to 11 levels of refinement (giving a solid angle of 1.5\xten{-6} sr.). For our robust goal-based metric, we use FP$_1$ as a surrogate, with a constant filter value of $\Sigma_\textrm{f}=1$. 
\begin{figure}[th]
\centering
\subfloat[][Error vs CDOFs]{\label{fig:2D_void_convg}\includegraphics[width =0.47\textwidth]{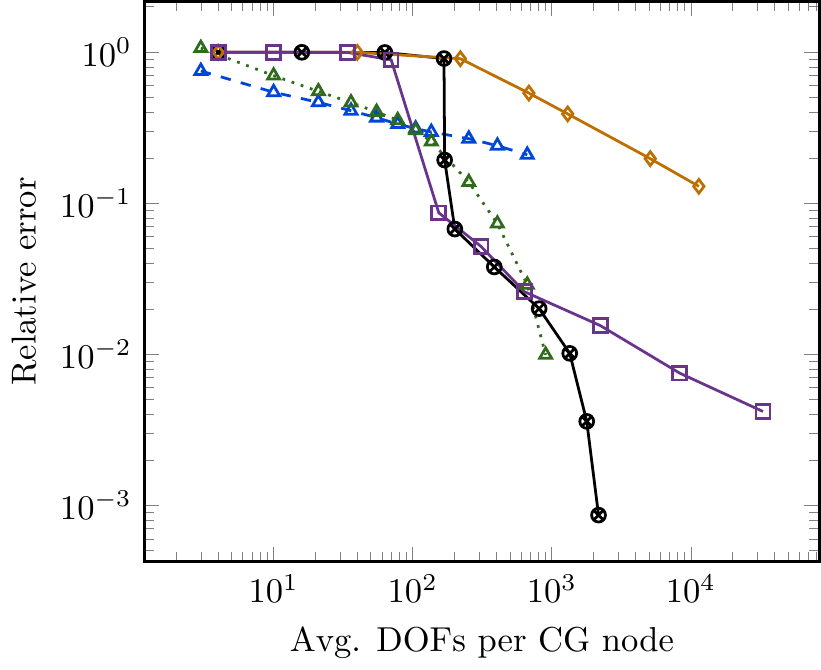}}
\subfloat[][Error vs total runtime]{\label{fig:2D_void_time}\includegraphics[width =0.47\textwidth]{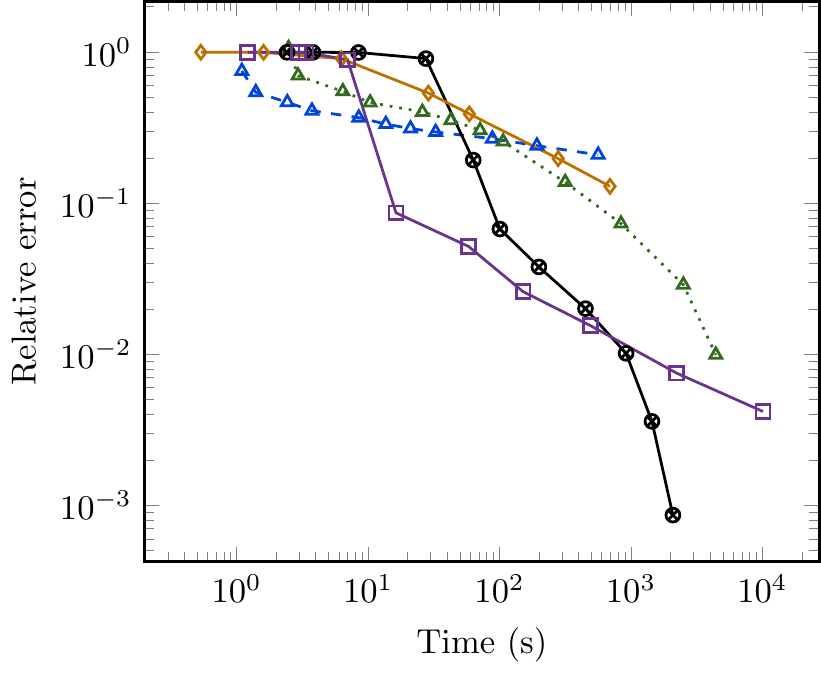}}
\caption{Comparison of the relative error of the detector response, for different angular discretisations for the 2D void problem. The \textcolor{gaylordpurple}{\Square} is non-standard Haar wavelets with fixed angular refinement between $\mu \in [0, 1]$ and $\omega \in [1.47976, 1.661832]$, the dashed \textcolor{matlabblue}{$\triangle$} is uniform FP$_n$ with $\Sigma_{\textrm{f}}=1$, the dotted \textcolor{foliagegreen}{$\triangle$} is uniform FP$_n$ with $\Sigma_{\textrm{f}}=0.1$, \textcolor{deludedorange}{$\diamond$} uniform LS P$^0$ FEM and the \textcolor{black}{$\otimes$} are goal-based adapted non-standard Haar wavelets with robust error target 1\xten{-3} and with reduced tolerance solves.}
\label{fig:2D_void_result}
\end{figure}

The goal for our error metric is the average flux in the detector region at the end of the duct. \fref{fig:2D_void_result} shows the results from our adaptive scheme on this problem when compared to several other angular discretisations. As expected, \fref{fig:2D_void_convg} shows that the relative error in our different NRI angular discretisations, namely the fixed refinement with Haars, the goal-based robust Haars and the uniform LS P$^0$ FEM are one with low levels of refinement, given ray-effects in this problem. The FP$_n$ method however, even at FP$_1$ is in the asymptotic regime given this problem, and with a uniform filter of $\Sigma_\textrm{f}=0.1$ converges well. This helps confirm that FP$_n$ discretisations can form a good surrogate to improve our error metric. We can also see that the uniform LS P$^0$ FEM performs poorly in this problem per DOF when compared with the fixed refinement Haars, again as would be expected in such a duct problem.
\begin{figure}[th]
\centering
\includegraphics[width =0.57\textwidth]{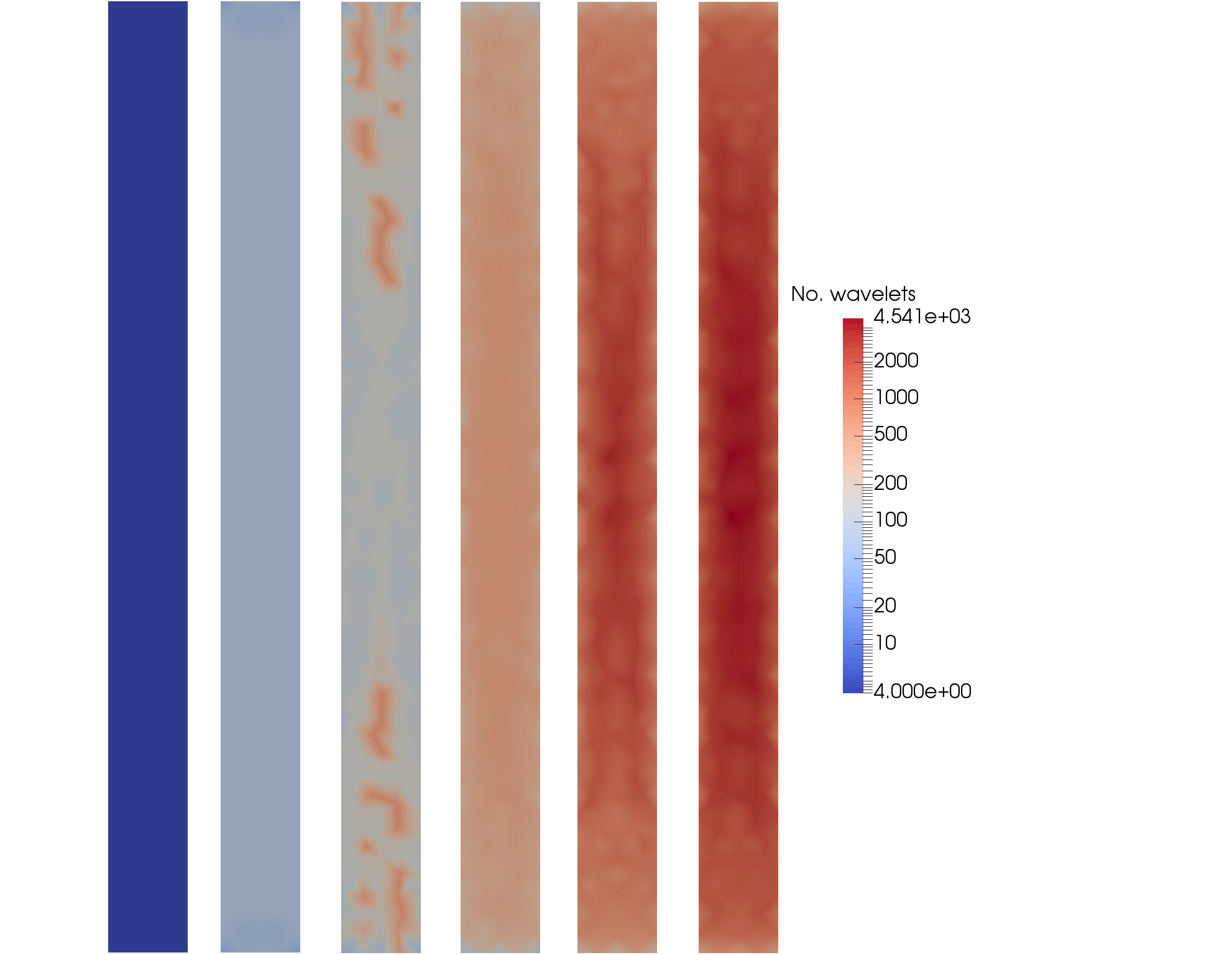}
\caption{Number of wavelets across the spatial domain for the 2D void problem, plotted on the CG mesh, on the different steps of the goal-based angular adaptivity with robust error target 1\xten{-3} (from \fref{fig:2D_void_result}). The 1st adapt step is shown on the left, moving to the right gives the 3rd, 5th, 7th, 9th and 11th adapt steps.}
\label{fig:2D_void_no_angles}
\end{figure}

Importantly we can see our adapted robust Haars refine correctly in this problem, with the error per adapt step almost matching the fixed refinement Haars. This is in contrast to a standard goal-based metric in this problem, which we cannot plot as it does not adapt at all in this problem. We also see that for higher number of adapt steps, the adapted robust Haars outperform the fixed refinement; this is because the adaptivity can focus resolution primarily in the centre of the duct as the particles are not particularly well collimated in such a short duct. \fref{fig:2D_void_no_angles} shows where in space the adaptivity has placed resolution and we can see that on the 11th adapt step, there is more resolution applied in the centre of the duct. This is in contrast to the fixed refinement which has applied the same angular resolution everywhere in space. 

\fref{fig:2D_void_time} shows the runtime of our method and like \fref{fig:2D_void_convg} we can see our adapted robust Haars performing well, producing a relative error of $\sim$1\xten{-3} roughly two orders of magnitude quicker than the fixed refinement Haars. Although the FP$_n$ method performed well per DOF and is competitive with the adapted robust Haars, \fref{fig:2D_void_time} shows that they are far more expensive to compute. Again this is to be expected given the cost of computing an FP$_n$ solution should scale like $\mathcal{O}(n^2)$ with angular refinement and helps highlight the importance of only using FP$_n$ as a surrogate solution to improve our error metric. In particular, computing our forward and adjoint FP$_1$ surrogates for the robust Haar adapt only takes 2 seconds in this problem, which is a small fraction of the total Haar adapt runtime. 

\begin{figure}[th]
\centering
\subfloat[][Forward]{\label{fig:2D_void_ratio}\includegraphics[width =0.3\textwidth]{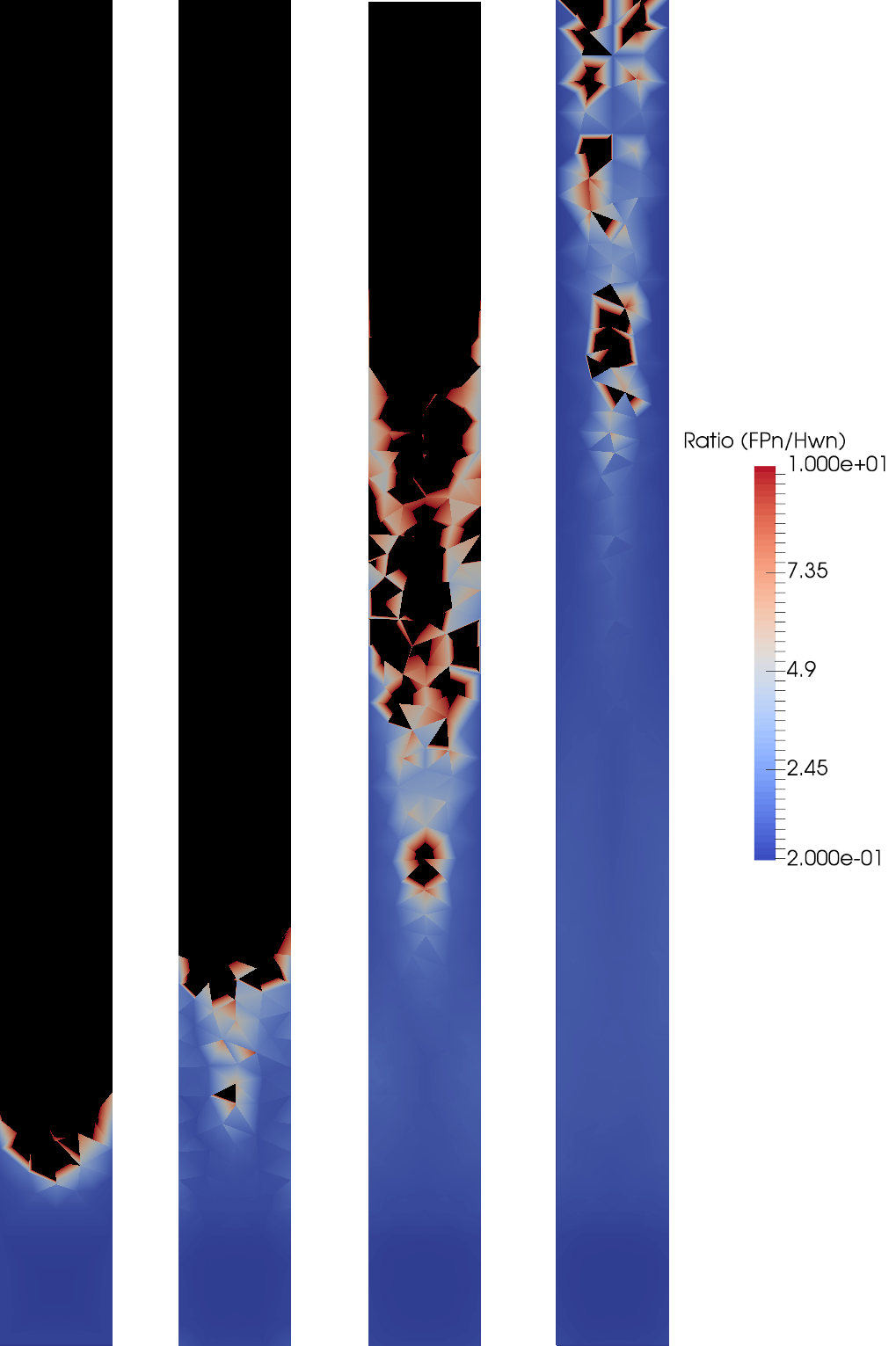}} \hspace{0.1cm}
\subfloat[][Adjoint]{\label{fig:2D_void_adjoint_ratio}\includegraphics[width =0.2\textwidth]{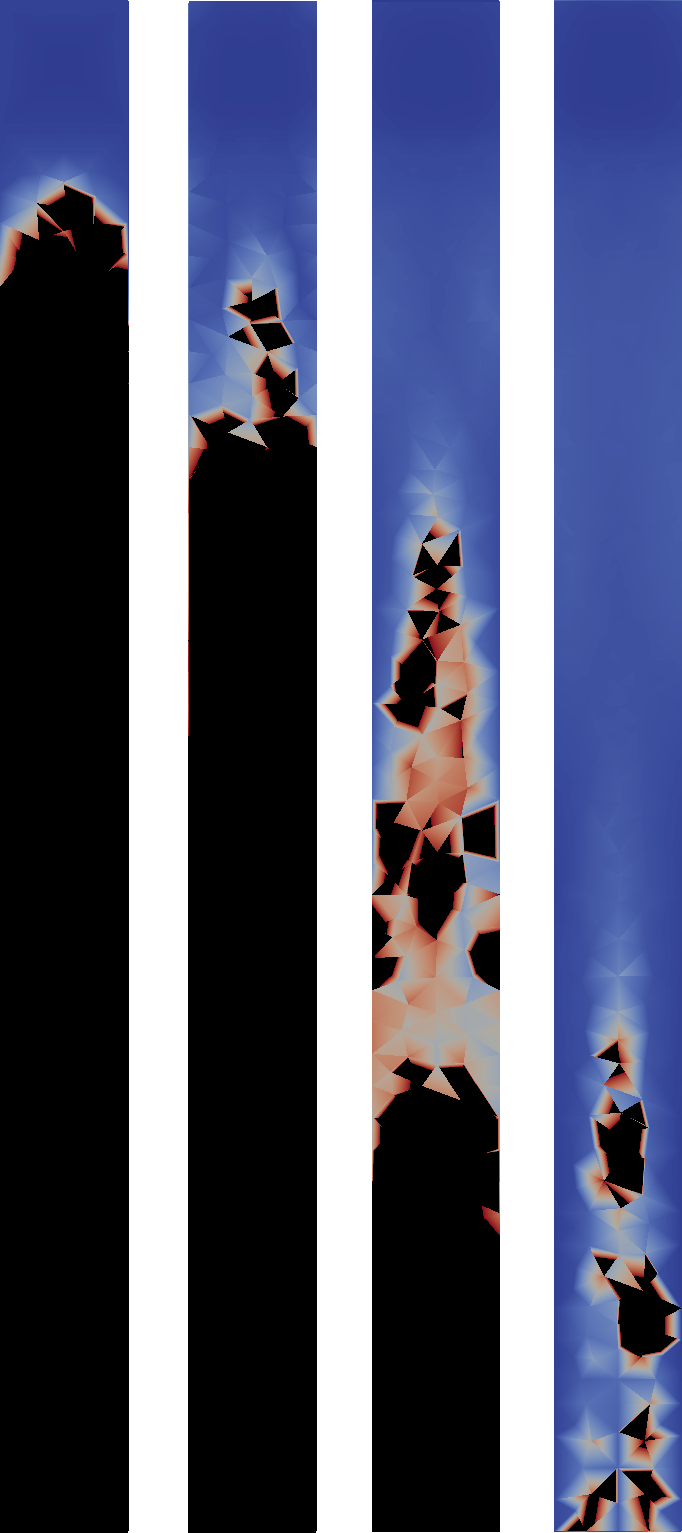}}
\caption{Ratio of the absolute value of the FP$_1$ scalar fluxes to the adapted Haar wavelet scalar fluxes in the 2D void problem, for the first four adapt steps (the first step is on the left). If the ratio is greater than 10, the nodes have been shaded black. These black nodes represet areas that are ``underresolved''.}
\label{fig:2D_void_ratios}
\end{figure}

To further confirm that our method for determining underresolved nodes is effective, \fref{fig:2D_void_ratios} shows the nodes in this problem that trigger using the mapped FP$_1$ solution in each of the individual forward and adjoint solutions, over the first four adapt steps shown in \fref{fig:2D_void_result}. We can see that our heuristic correctly captures the areas of the spatial domain affected by ray-effects and that as our adaptivity progresses, the number of nodes considered underresolved shrinks. Algorithm \ref{alg:one} considers both the forward and adjoint solutions when labelling a node as underresolved, and the percentage of nodes labelled as such are 100\%, 100\%, 71\% and 12\% for the first four adapt steps shown in \fref{fig:2D_void_ratios}. Beyond this step, 0\% of nodes are labelled underresolved and our robust error metric reverts to the non-robust, as we are in the asymptotic regime and the non-robust metric suffices.
\begin{table}[ht]
\centering
\begin{tabular}{ l c c c c c c c c c c c c}
\toprule
\textbf{Adapt step:} & \textbf{1} & \textbf{2} & \textbf{3} & \textbf{4} & \textbf{5} & \textbf{6} & \textbf{7} & \textbf{8} & \textbf{9} & \textbf{10} & \textbf{11}\\
\midrule  
Non-robust metric & 5.8\xten{-14} & - & - & - & - & - & - & - & - & - & - \\
Robust metric & 1.03 & 1.25 & 0.84 & 0.16 & 0.74 & 3.22 & 9.87 & 68.92 & 13.63 & 9.15 & 8.18\\
\bottomrule  
\end{tabular}
\caption{Effectivity index for the goal-based adapted discretisation shown in \fref{fig:2D_void_result}, for the 2D void problem.}
\label{tab:2D_void_effec}
\end{table}

Furthermore, \tref{tab:2D_void_effec} shows the effectivity index for the robust Haar adapt and we can see that it is non-zero for all the adapt steps. In the pre-asymptotic regime when the FP$_n$ solution is used throughout most of the domain, we see the effectivity index is very close to 1, while the non-robust metric is (machine) zero and adaptivity does not occur. This confirms that using our FP$_n$ solution to bootstrap our error metric when ray-effects are present is effective. The effectivity index in the asymptotic regime becomes much larger as it has reverted to using the simple diagonal matvec described by \eref{eq:disc_resid_subgrid} and discussed in \cite{Dargaville2019}, but as \fref{fig:2D_void_result} shows we are still adapting in the correct regions. 

\begin{figure}[th]
\centering
\subfloat[][Forward FP$_{1}$ solution with $\Sigma_\textrm{f}=1$.]{\label{fig:fpn_forward_middle_step_11}\includegraphics[width =0.32\textwidth]{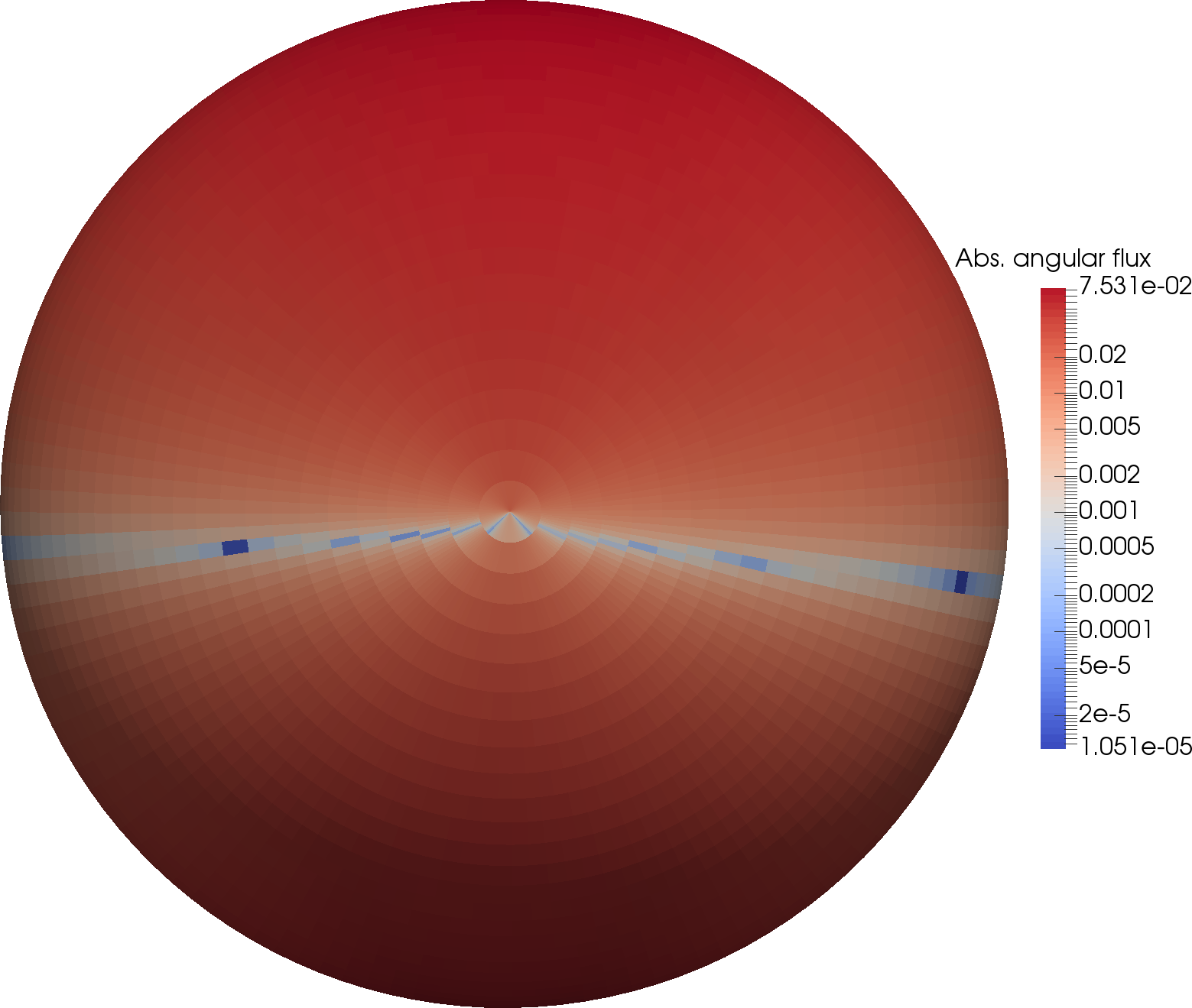}} \hspace{0.1cm}
\subfloat[][Forward solution with adapted wavelets after 3 adapt steps with robust error target 1\xten{-3}.]{\label{fig:forward_middle_step_3}\includegraphics[width =0.32\textwidth]{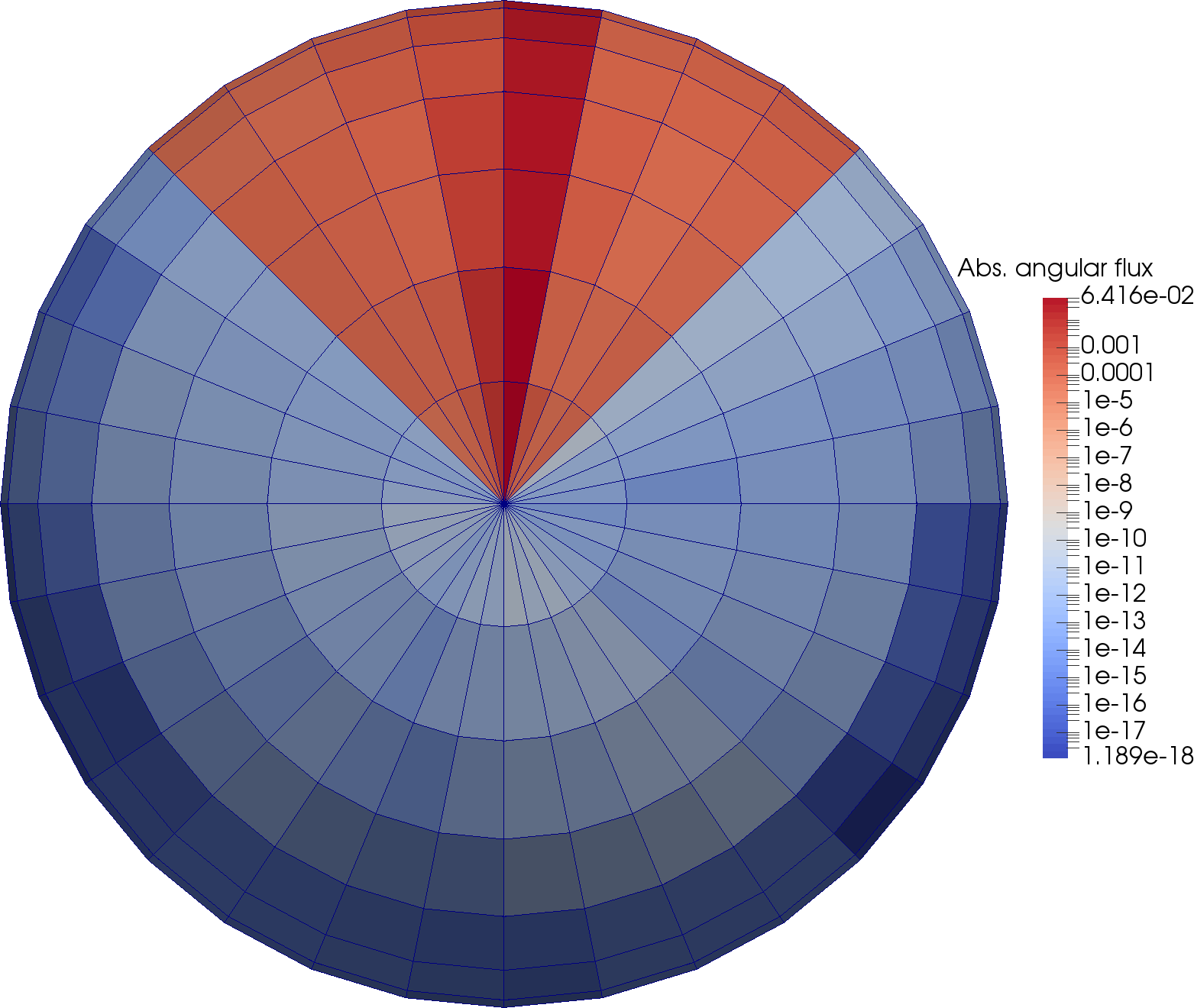}} \hspace{0.1cm}
\subfloat[][Forward solution with adapted wavelets after 11 adapt steps with robust error target 1\xten{-3}.]{\label{fig:forward_middle_step_11}\includegraphics[width =0.3\textwidth]{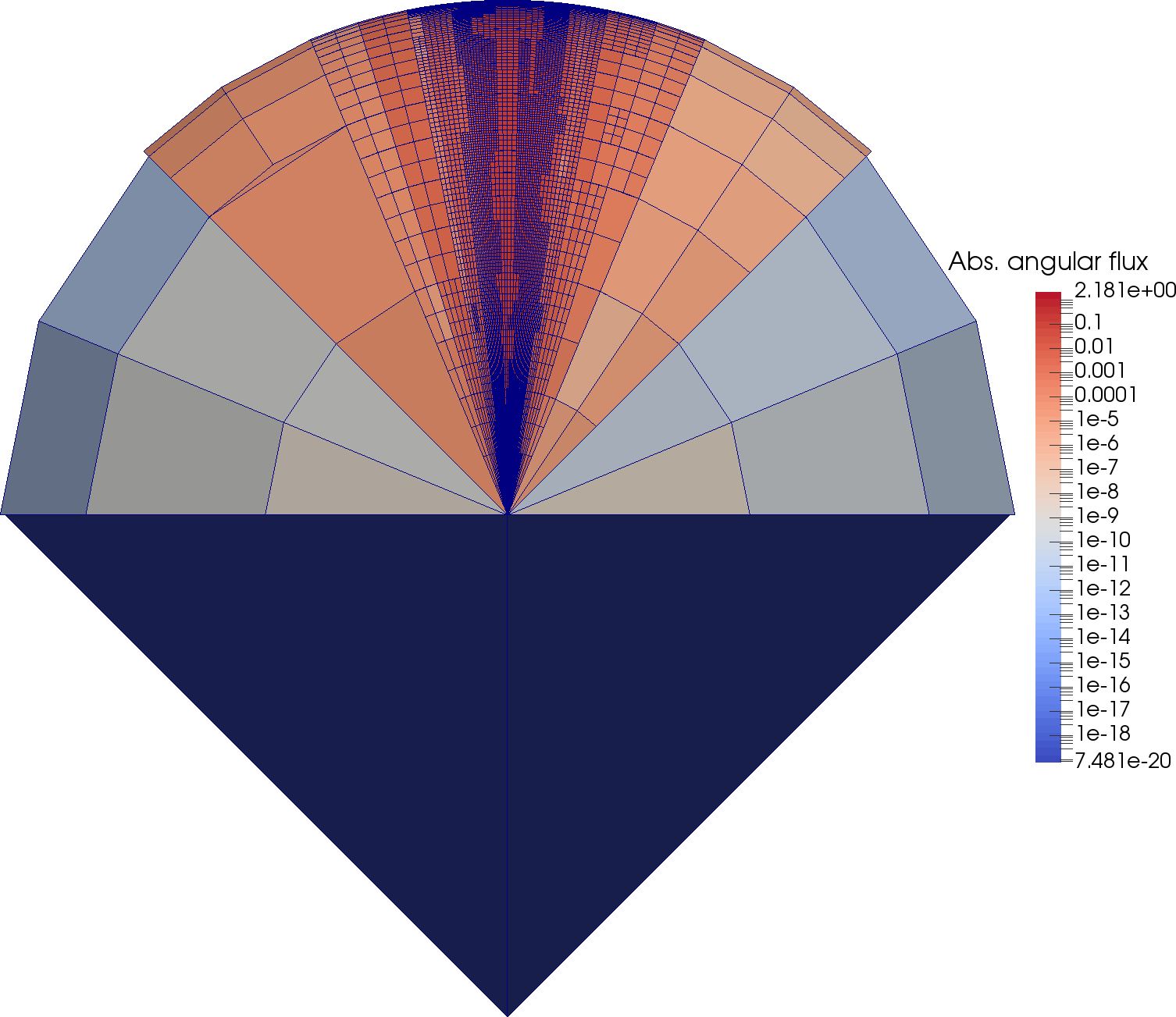}}
\caption{Absolute value of the angular flux in the 2D void problem at spatial position $x=0.390625$, $y=6.25$ (from \fref{fig:2D_void_result}). All angular discretisations are on the $r=1$ sphere, but have been projected onto faceted polyhedra for ease of visualisation. The camera is pointed in the $-z$ direction.}
\label{fig:ang_flux}
\end{figure}

We now plot the absolute value of the angular flux from our FP$_1$ surrogate and our adapted robust Haars close to the midpoint of the duct in \fref{fig:ang_flux}. We can see that the FP$_1$ solution in \fref{fig:fpn_forward_middle_step_11} appears to be almost uniform on the sphere, but highest flux is correctly pointed down the duct. Interestingly, we can still see oscillations in the solution, and the solution is in fact negative in the $-y$ direction (the two small regions of $\sim 1$\xten{-8} correspond to where the solution turns negative). \fref{fig:forward_middle_step_3} shows the angular flux of our adapted robust Haars in the third adapt step and \fref{fig:2D_void_ratios} shows that the middle of the duct is still considered underresolved at this point. This means that the error metric is forcing the wavelets to adapt to the FP$_1$ solution shown in \fref{fig:fpn_forward_middle_step_11}. This is why the ``adapted'' discretisation shown in \fref{fig:forward_middle_step_3} is uniform; the lack of anisotropy in the angular flux at FP$_1$ and the oscillations are triggering uniform Haar refinement on the sphere. We can see this in \fref{fig:2D_void_convg}, where the adapted robust Haars for the first three adapt steps are using more DOFs than the fixed refinement Haars, indicating they are not refining anisotropically on the sphere. 

This is a consequence of using such a low order FP$_n$ solution as our surrogate. This is an important point, as using FP$_1$ in this problem is triggering the same uniform refinement that using a diffusion solution would as a surrogate. We could have carefully tuned our thresholding parameter $\tau$ in our error metric to focus on the slightly higher angular flux pointing down the duct, and produce an anisotropic adapt in this problem. We instead chose to highlight this point and discuss this further in \secref{sec:2D void problem - 100}. Finally \fref{fig:forward_middle_step_11} shows that once we are in the asymptotic regime on adapt step 11, where 0\% of the nodes are marked as underresolved, the anisotropic adaptivity has resumed and we see heavy refinement down the direction of the duct (and mainly around the polar region, which is important given the symmetry in this 2D problem). In particular, we can see the $-y$ direction has coarsened the unnecessary uniform refinement triggered by the FP$_1$ surrogate shown in \fref{fig:forward_middle_step_3}; this coarsening is visible in step 5 of the adaptive robust Haars in \fref{fig:2D_void_convg}, as we increase the level of refinement (and hence decrease the relative error) the net NDOFs goes down. To further investigate this behaviour, we turn to our next example problem. 
\subsection{2D void problem - 100}
\label{sec:2D void problem - 100}
\begin{figure}[th]
\centering
\subfloat[][Relative error of the detector response vs CDOFs]{\label{fig:2D_void_100_convg}\includegraphics[width =0.4\textwidth]{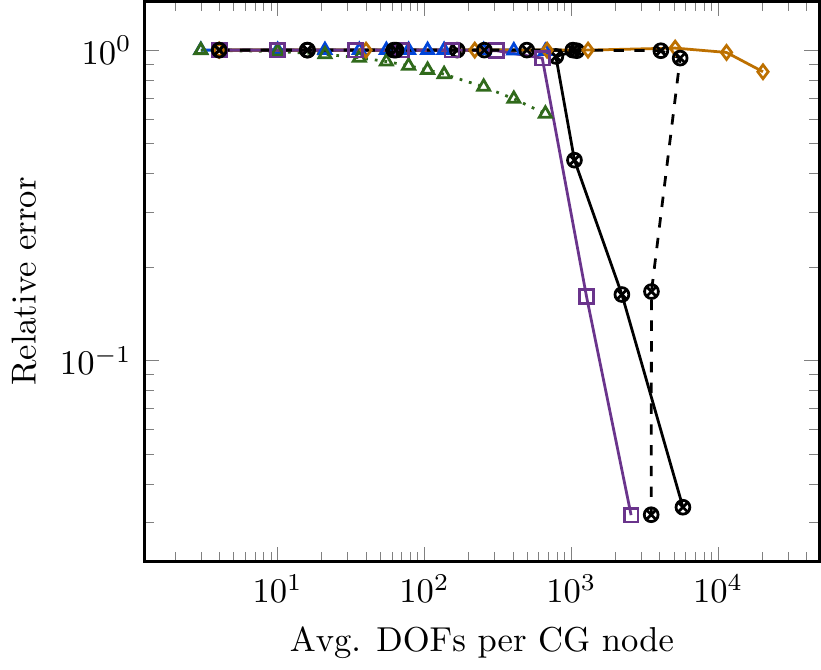}}
\subfloat[][Relative error of the detector response vs total runtime]{\label{fig:2D_void_100_time}\includegraphics[width =0.4\textwidth]{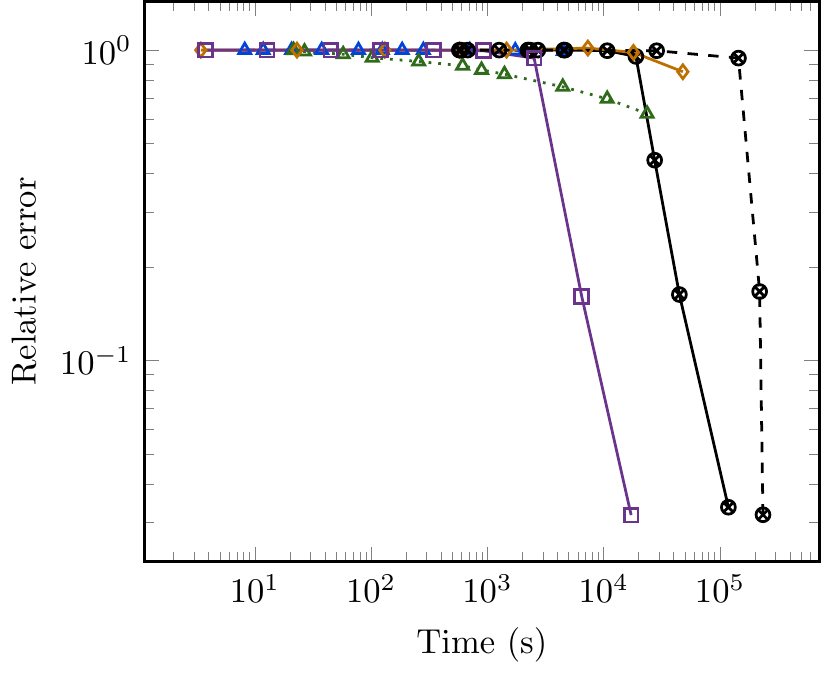}}\\
\subfloat[][Absolute value of the average flux in the detector vs CDOFs]{\label{fig:2D_void_100_response}\includegraphics[width =0.4\textwidth]{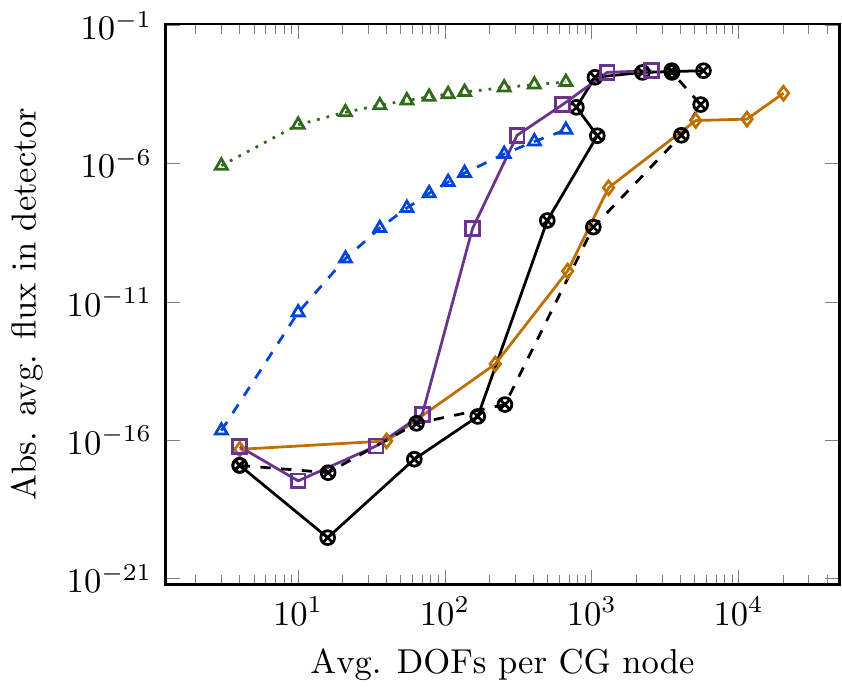}}
\caption{Comparison of the results for different angular discretisations for the 2D void problem of length 100. The \textcolor{gaylordpurple}{\Square} is non-standard Haar wavelets with fixed angular refinement between $\mu \in [0, 1]$ and $\omega \in [1.561, 1.5807]$, the dashed \textcolor{matlabblue}{$\triangle$} is uniform FP$_n$ with $\Sigma_{\textrm{f}}=1$, the dotted \textcolor{foliagegreen}{$\triangle$} is uniform FP$_n$ with $\Sigma_{\textrm{f}}=0.1$, \textcolor{deludedorange}{$\diamond$} uniform LS P$^0$ FEM and the \textcolor{black}{$\otimes$} are goal-based adapted non-standard Haar wavelets with robust error target 1\xten{-6} and with reduced tolerance solves. The dashed \textcolor{black}{$\otimes$} uses FP$_9$ with $\Sigma_\textrm{f}=0.1$ as a surrogate.}
\label{fig:2D_void_100_result}
\end{figure}
In this example, we extend the length of the duct to 100cm, with the same source and detector and pure vacuum of the problem in \secref{sec:2D void problem} (where now the total length of the domain is 102cm), discretised with a 4334 element unstructured triangular mesh (2580 CG nodes and 13,002 DG nodes). Again we compare against \textit{a priori} fixed refinement between $\mu \in [0, 1]$ and $\omega \in [1.561, 1.5807]$ with our non-standard Haar wavelets. The reference solution used is the non-standard fixed refinement Haars up to 11 levels of refinement (giving a solid angle of 1.5\xten{-6} sr.). As in the previous example, using the non-robust goal-based error metric does not trigger refinement in this problem. For our robust goal-based metric, we use FP$_{21}$ as a surrogate, with a constant filter value of $\Sigma_\textrm{f}=1$ and an extra adapt step at max order.

To begin, we can see in \fref{fig:2D_void_100_result} that the relative error in this problem for the fixed refinement Haars, the LS P$^0$ FEM and the adapted robust Haars are one for considerably more refinement steps than when the duct was of length 10, as would be expected. It takes a much smaller solid angle to get a detector response in this problem; \fref{fig:2D_void_100_response} shows it takes approximately $4\xtenm{2}$ DOFs to record any significant response for the NRI angular discretisations. Similarly we can see that the fixed refinement Haars and the adapted robust Haars outperform the uniformly refined LS P$^0$ both per DOF and in total runtime, shown in Figures \ref{fig:2D_void_100_convg} \& \ref{fig:2D_void_100_time}, respectively. 
\begin{figure}[th]
\centering
\subfloat[][Forward FP$_{21}$ solution with $\Sigma_\textrm{f}=1$.]{\label{fig:void_100_fp21_forward}\includegraphics[width =0.32\textwidth]{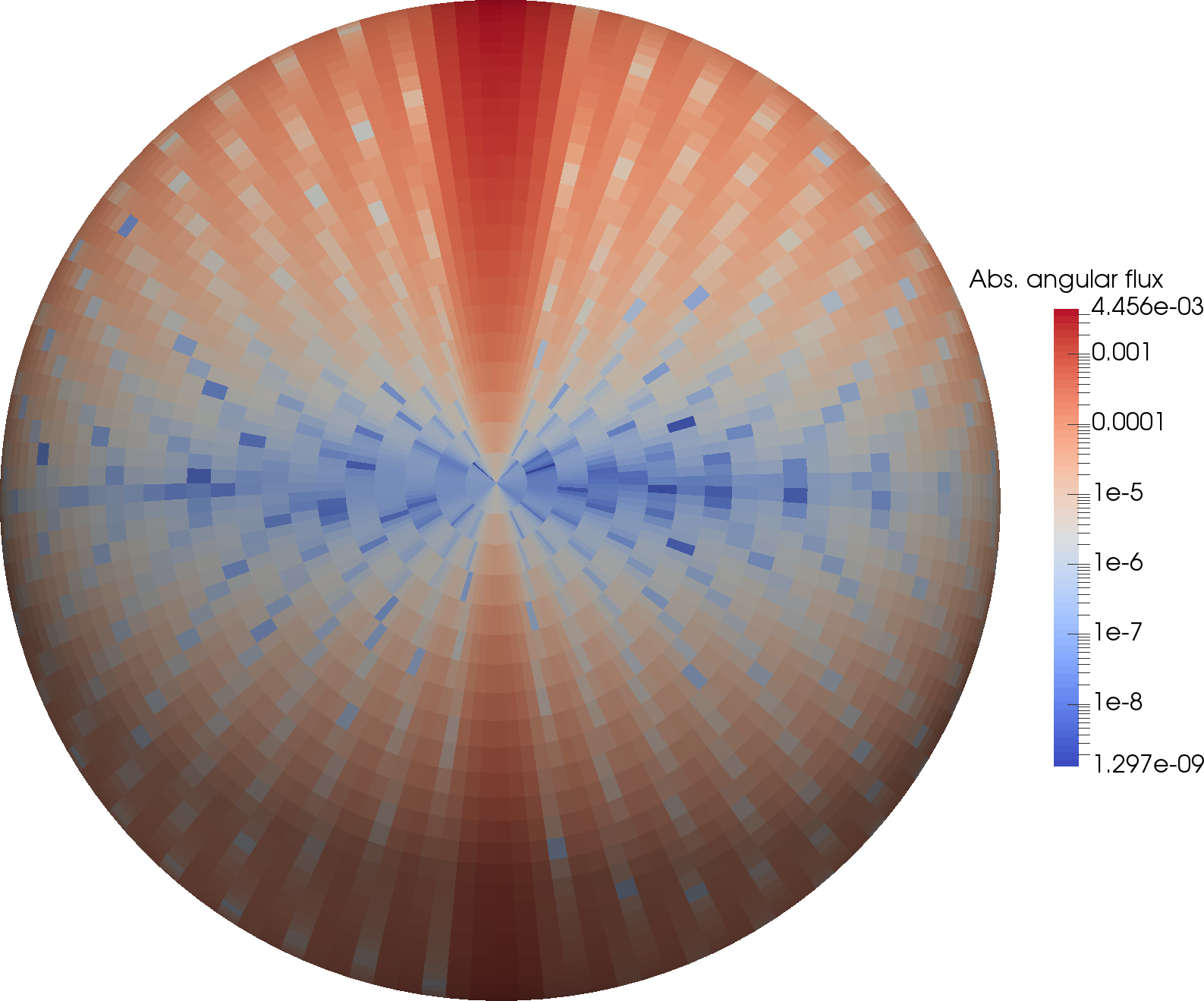}} \hspace{0.1cm}
\subfloat[][Forward solution with adapted wavelets after 5 adapt steps with robust error target 1\xten{-6}. The element boundaries have been coloured white for visiblity.]{\label{fig:void_100_forward_5}\includegraphics[width =0.32\textwidth]{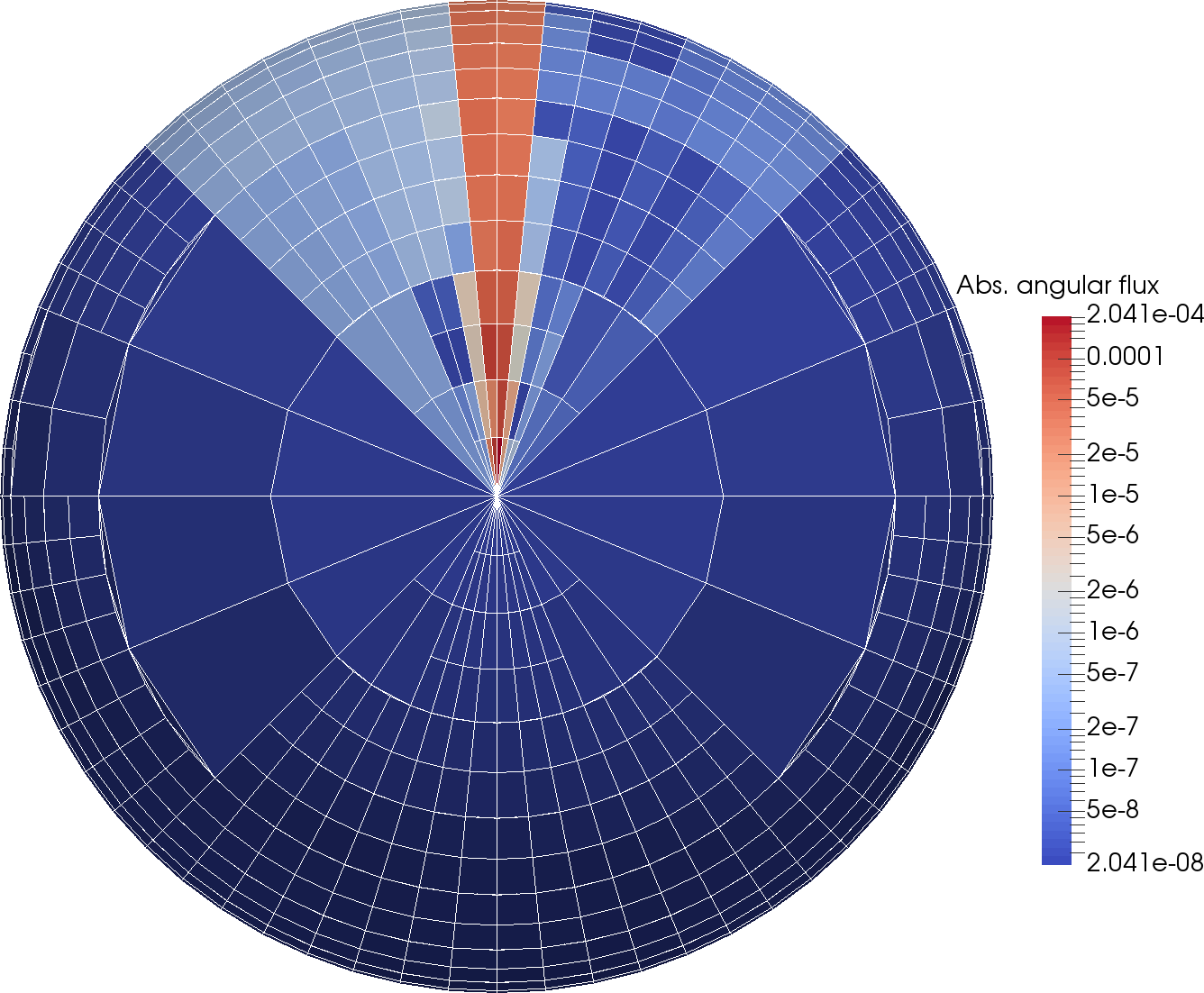}} \hspace{0.1cm}
\subfloat[][Forward solution with adapted wavelets after 10 adapt steps with robust error target 1\xten{-6}.]{\label{fig:void_100_forward_10}\includegraphics[width =0.3\textwidth]{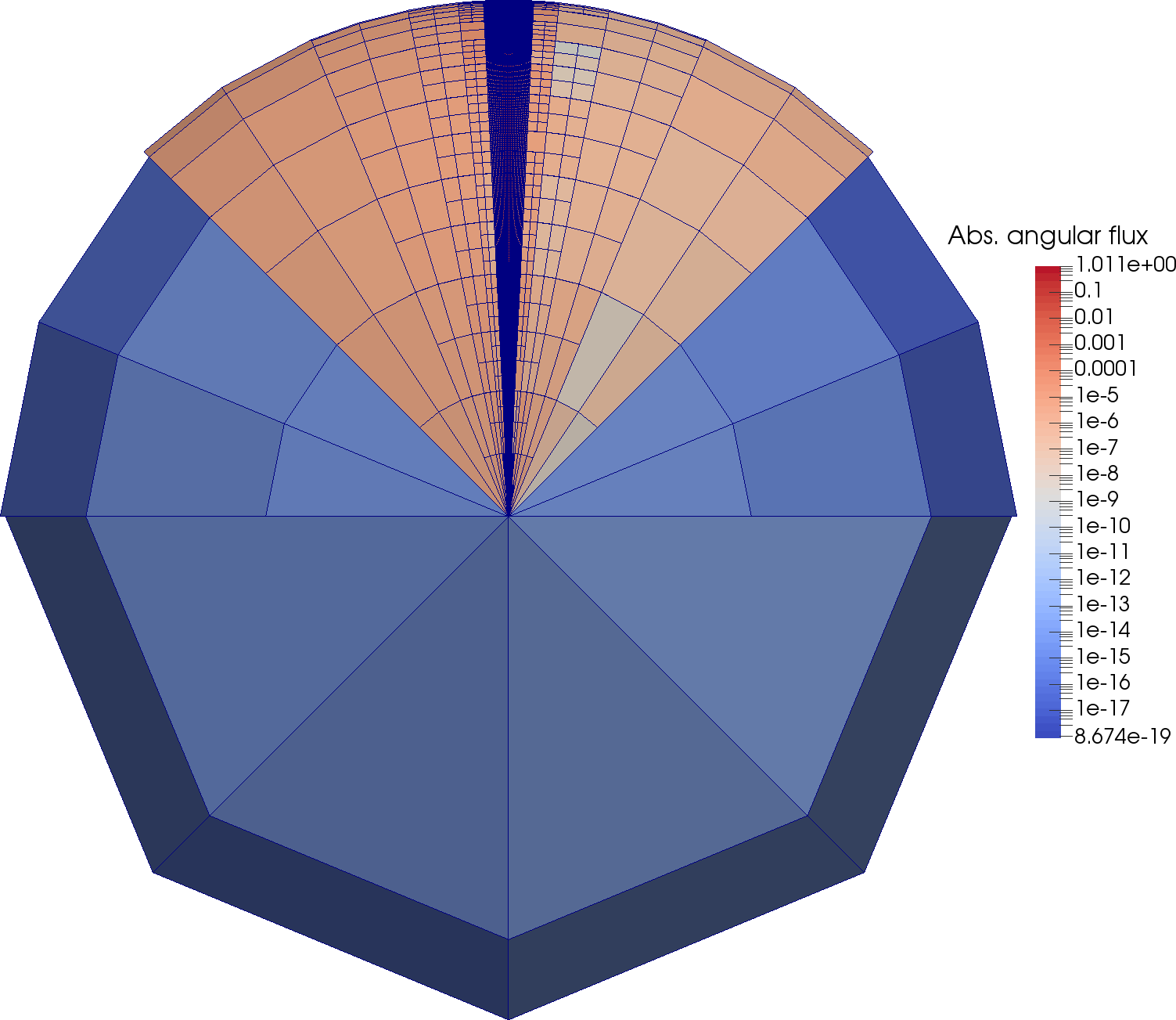}}
\caption{Absolute value of the angular flux in the 2D void problem of length 100 at spatial position $x=0.608$, $y=51$ (from \fref{fig:2D_void_100_result}). All angular discretisations are on the $r=1$ sphere, but have been projected onto faceted polyhedra for ease of visualisation. The camera is pointed in the $-z$ direction.}
\label{fig:ang_flux_100}
\end{figure}

\fref{fig:2D_void_100_response} shows that the FP$_n$ records a detector response far before the NRI discretisations in this problem, even at low order and with a heavy filter. Importantly we can see that it is better to run with a smaller filter at lower order than use a stronger filter and high order discretisation, as the FP$_1$ discretisation with $\Sigma_\textrm{f}=0.1$ records the same detector response as FP$_{21}$ with $\Sigma_\textrm{f}=1$. This highlights the importance of not overfiltering, as discussed in \cite{Dargaville2019a}. Of course using a higher filter value decreases the condition number and reduces the iteration count in the iterative method, so we must strike a balance in these problems. Importantly however, for both of the FP$_n$ results shown, they are in the asymptotic regime far earlier than the NRI angular discretisations. \fref{fig:2D_void_100_time} also helps emphasise that increasing the order of the FP$_n$ solutions is not a practical solution method in these streaming problems given the $\mathcal{O}(n^2)$ cost. We can see in \fref{fig:2D_void_100_time} that once the fixed refinement Haars and the robust adapt Haars record a detector response, they are far more efficient to compute than either FP$_n$ discretisation in terms of converging the detector response.

Clearly however, the poorly converged FP$_{21}$ surrogate used by the adapted robust Haars is sufficient to bootstrap the NRI adaptivity and causes the Haars to adapt correctly; we see in \fref{fig:2D_void_100_convg} that the adapted Haars match the error of the fixed refinement Haars at each level of refinement. The adapted Haars use slightly more DOFs in the asymptotic region and takes a fixed amount of extra time per refinement step compared to the fixed refinement. We should also note that computing our forward and adjoint FP$_{21}$ surrogates only cost an additional $\sim\,1$\% in runtime for our robust adapted Haars with 9 adapt steps. 

We would not expect the adapted Haars to outperform the fixed refinement in this problem per DOF, as the adaptivity cannot remove extra unimportant DOFs at the ``edges'' of the duct compared to the fixed refinement as in \secref{sec:2D void problem}, given the greater distance between source and detector. We also do not expect the adapted Haars to be quicker in this problem, given we are only required to do one linear solve to with the fixed refinement, as we know where to refine \textit{a priori}, compared with 23 linear solves for the adapted robust Haars with an extra step and max refinement level of 10. 

\fref{fig:2D_void_100_result} also shows the results from where we have used a different surrogate solution, namely FP$_9$ with a smaller constant filter value of $\Sigma_\textrm{f}=0.1$. We can see in the pre-asymptotic regime in \fref{fig:2D_void_100_response} that this has used more DOFs compared to the FP$_{21}$ surrogate, but that in the asymptotic regime the error and NDOFs matches the FP$_{21}$ surrogate well. 

One interesting point to note in \fref{fig:2D_void_100_convg} is that the FP$_9$ robust Haar adapt shows the the number of DOFs decreases, while the error also decreases, as the refinement level is increased. We can see in \fref{fig:2D_void_100_response} however that both robust adapted Haar solutions decrease their NDOFs once they reach their 7th adapt step. This is the transition point between the pre-asymptotic regime where the FP$_n$ surrogate solution is being used in the error metric, and the asymptotic regime where the Haar wavelets have sufficiently resolved the detector. The decrease in NDOFs when this occurs implies that unnecessary resolution has been applied in the pre-asymptotic regime, and that this can be coarsened in subsequent adapt steps. 

This is indeed the case and to examine this further, \fref{fig:ang_flux_100} plots the angular flux for the FP$_n$ and Haar solutions at the mid-point of the duct. We can see the forward FP$_{21}$ surrogate solution in \fref{fig:void_100_fp21_forward}. We used a higher FP$_n$ order as surrogate in this solution when compared to \secref{sec:2D void problem} so that the surrogate would not prompt uniform refinement and we can correctly see a strong peak pointing down the direction of the duct, but oscillations in the solution are also evident. The filtering of spherical harmonics does not completely prevent Gibbs-like oscillations in the solution, it merely dampens them sufficiently to give a (near) constant condition number. In particular, we see a large oscillation pointing in the $-y$ direction. We should note this oscillation is spurious given \fref{fig:void_100_fp21_forward} shows the forward solution; a similar oscillation in the opposite (wrong) direction is visible in the adjoint solution. \fref{fig:void_100_forward_5} shows that in the pre-asymptotic regime for the Haars after 5 adapt steps, given the tolerance of 1\xten{-6} used in the adapt process, the Haar solution has adapted to capture both the correct $+y$ peak in the flux, but also the spurious $-y$ oscillation. We could have adjusted the adapt tolerance to try and hide this effect (as the oscillation is smaller than the true peak), or have further filtered the FP$_n$ solution to smooth out these oscillations, but as we have chosen a sufficiently high resolution FP$_n$ surrogate (and sufficiently ``small'' filter value), adapting to the spurious oscillation in the pre-asymptotic region only costs us roughly three times the NDOFs when compared to the fixed refinement Haars, which uses fewer DOFs than uniform refinement.

As such this is better than the uniform refinement seen in \secref{sec:2D void problem} and we can see that this excess resolution is dropped once the Haars reach the asymptotic regime; \fref{fig:void_100_forward_10} shows the robust Haar forward solution on the 10th adapt step and we can see very little resolution in the $-y$ direction. This highlights that we simply have to pick an FP$_n$ surrogate solution and adapt tolerance that is ``good enough'' to prevent uniform refinement to get good performance from our robust metrics. If we do not, for example when we used the FP$_9$ surrogate with $\Sigma_\textrm{f}=0.1$, then we pay a cost in NDOFs and hence runtime for early adapt steps, but as \fref{fig:2D_void_100_result} shows our adapt process still produces a robust, anisotropically adapted solution in asymptotic region.
\begin{table}[ht]
\centering
\begin{tabular}{ l c c c c c c c c c c}
\toprule
\textbf{Adapt step:} & \textbf{1} & \textbf{2} & \textbf{3} & \textbf{4} & \textbf{5} & \textbf{6} & \textbf{7} & \textbf{8} & \textbf{9} & \textbf{10}\\
\midrule  
Non-robust metric & 5.06\xten{-14} & - & - & - & - & - & - & - & - & -\\
Robust metric, FP$_{21}$, $\Sigma_\textrm{f}=1$ & 0.05 & 0.12 & 0.24 & 0.45 & 0.38 & 0.33 & 8.05 & 66.9 & 268.4 & 13559.2 \\
Robust metric, FP$_{9}$, $\Sigma_\textrm{f}=0.1$ & 25.2 & 48.9 & 118.4 & 176.8 & 189.4 & 56.3 & 20.1 & 227.98 & 1599.3 & - \\
\bottomrule  
\end{tabular}
\caption{Effectivity index for the goal-based adapted discretisation shown in \fref{fig:2D_void_100_result}, for the 2D void problem of length 100.}
\label{tab:2D_void_100_effec}
\end{table}

\tref{tab:2D_void_100_effec} shows the effectivity index for both the robust Haar adapts in \fref{fig:2D_void_100_result}, and similar to the results in \secref{sec:2D void problem} we see regardless of our choice of surrogate solution, the index is non-zero in the pre-asymptotic regime and increase in the asymptotic, unlike the non-robust metric. Given we have established our method is robust in pure vacuum problems with length/width ratios up to 1/100, we now move onto a problem that involves scattering and multiple streaming paths to further test the robustness of our method. 
\subsection{3D scatter box problem}
\label{sec:3D scatter box}
\begin{figure}[th]
\centering
\includegraphics[width =0.5\textwidth]{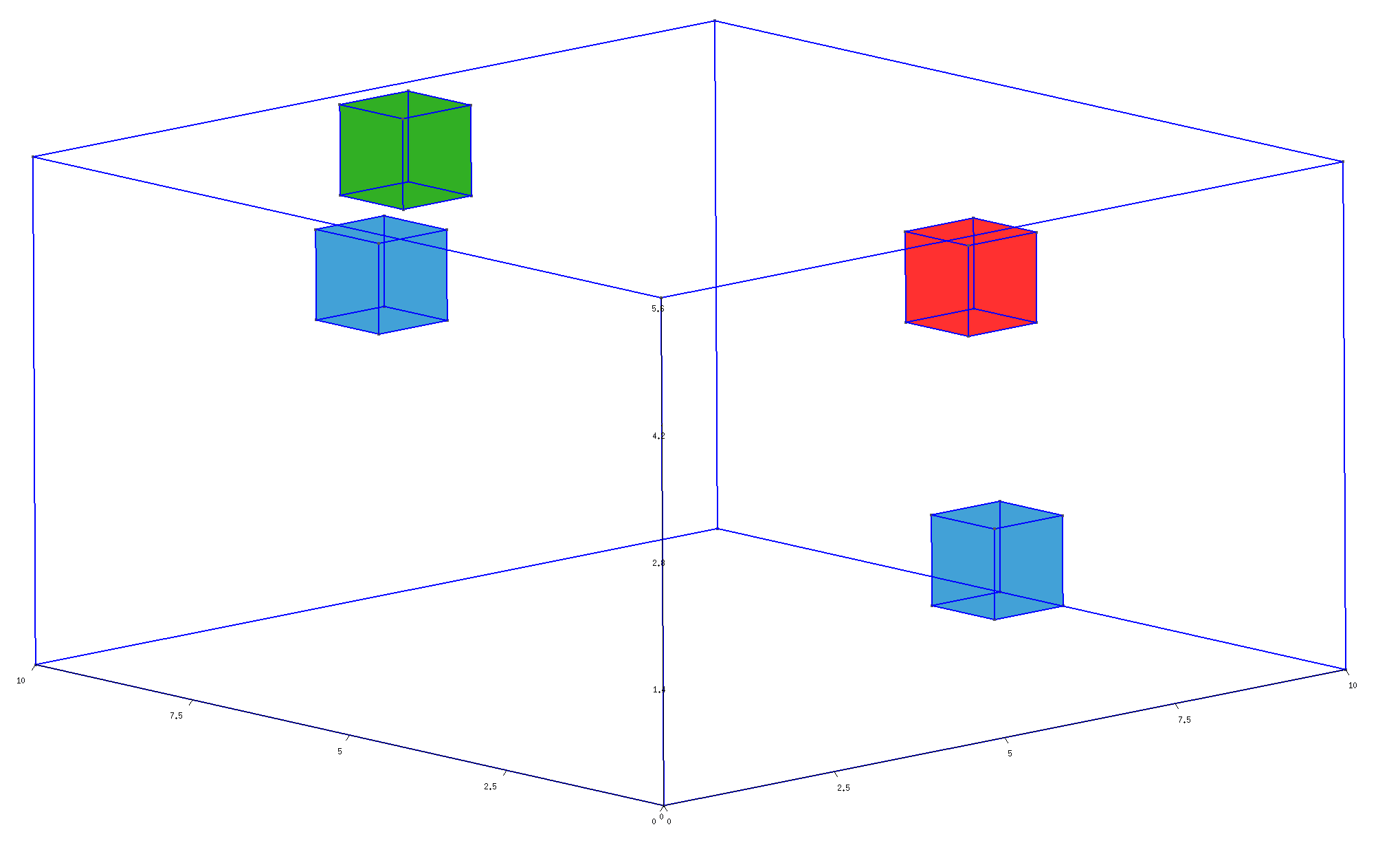}
\caption{Schematic of the 3D scatter box problem, with dimensions 10 x 10 x 6 cm. The red region is a 1 x 1 cube centred at (5, 0.5, 5) and is a source of strength 1, the green region is a 1 x 1 cube centred at (5, 9.5, 5) with the two blue cubes centred at (9.5, 5, 0.5) and (0.5, 5, 5) as pure scattering material (1 cm$^{-1}$). All regions except the blue scattering boxes are pure vacuum.}
\label{fig:3D_scatter_box_schem}
\end{figure}
The final example we test our adaptivity algorithm on is a 3D streaming problem with two small scattering regions, a schematic of which is shown in \fref{fig:3D_scatter_box_schem}. This example features five distinct paths through which particles can move from the source to the detector; the strongest of which is the direct path. The previous examples shown above highlight that our adaptivity algorithm handles the direct streaming path well, but the scattering boxes in this problem are designed to test the robustness of the method to other small contributions to the detector. We have also placed one of the scattering boxes such that the streaming path from source to box and then from box to detector is ``seen'' with our coarsest angular discretisation, H$_1$ (similar to S$_2$). As such we expect the non-robust error metric to trigger refinement, but only along that one path. The robust metric described in this work however should resolve all the streaming paths in this problem. 

We discretise this problem with an unstructured tetrahedral mesh with 21,650 elements (4886 CG nodes and 86,600 DG nodes). The fixed refinement non-standard Haar solution is \textit{a priori} refined between $\mu \in [-1, 1]$ and $\omega \in [0, 3.15]$. This is a large region in angle, but given the different streaming paths we cannot be more selective \textit{a priori} in this problem. Given this, we could not run high enough refinement level with the fixed refinement Haar to provide a sufficient reference to test our adaptive methods; as such we resort to using our robust adapted Haars as a reference in this problem; we use 11 adapt steps, with a maximum refinement level of 10 and a tolerance of 1\xten{-4} as our reference. This is at least three refinement levels and one order of magnitude smaller tolerance than any adapted solution we compare against and we verified that the solution produced matched the fixed refinement solution up to the highest level of refinement we could run which was 5 adapt steps. For the LS P$^0$ solution, we used a reference with 5328 elements in angle (similar to S$_{72}$). For our robust goal-based metric, we use the FP$_3$ as a surrogate with a constant filter value of $\Sigma_\textrm{f}=1$. 
\begin{figure}[th]
\centering
\includegraphics[width =0.5\textwidth]{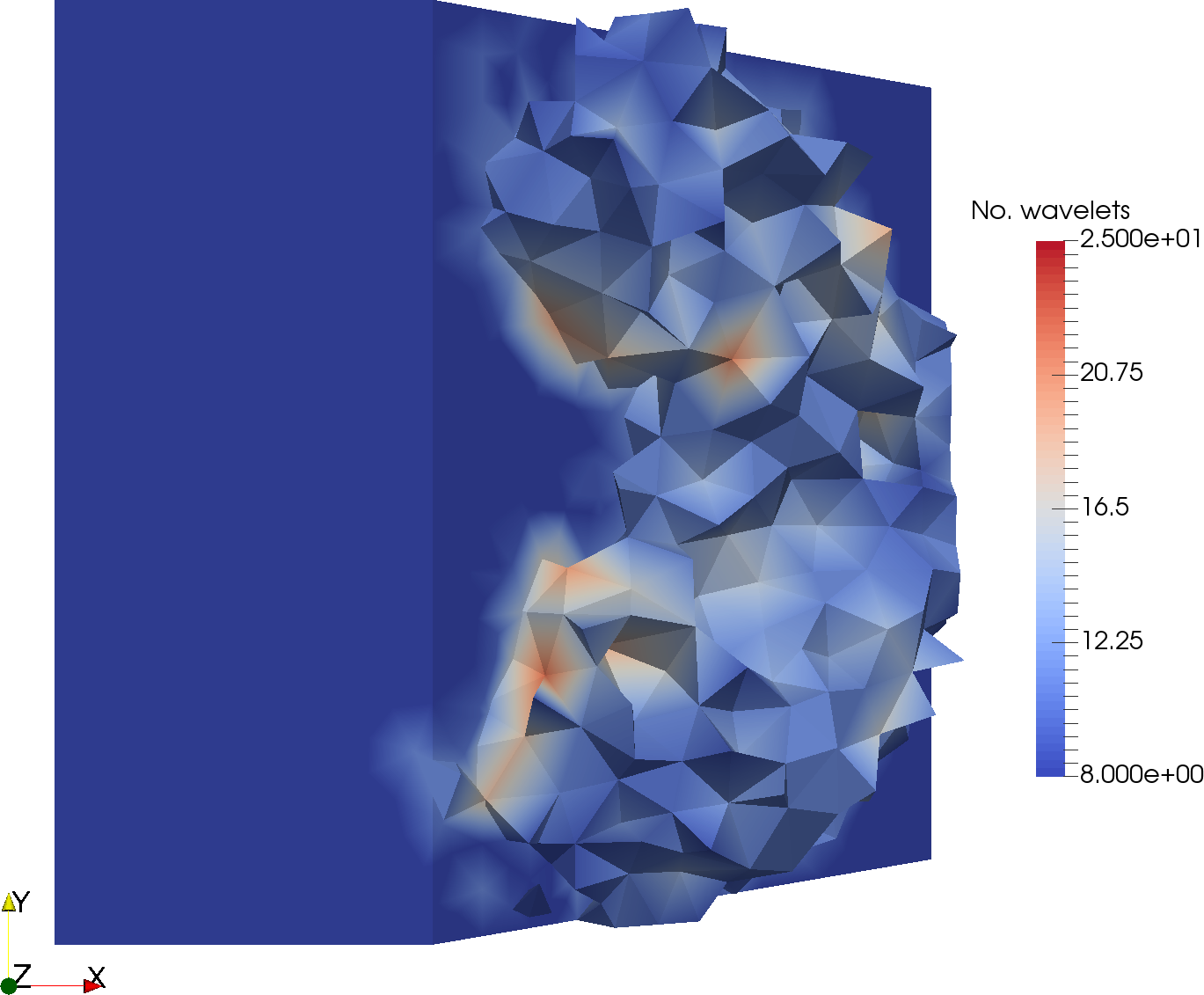}
\caption{Number of wavelets across the spatial domain for the 3D void problem, plotted on the CG mesh on the 2nd step of the goal-based (non-robust) angular adaptivity with error target 1\xten{-2} (the \textcolor{darksalmon}{x} line in \fref{fig:2D_void_100_result}). A cut and isosurface have been made in the visualisation, to show the only region where angular adaptivity has triggered refinement, between the source, the scattering region centred at (9.5, 5, 0.5) and the detector. The camera is pointed in the -z direction.}
\label{fig:not_robust_adapt}
\end{figure}

To begin, \fref{fig:not_robust_adapt} shows a visualisation of where the non-robust error metric has placed angular resolution across the domain after the second adapt step. We can see that as expected, it has refined only in the streaming path between the source, one of the scatter boxes and the detector. As mentioned this is because the H$_1$ discretisation aligns directly with this single streaming path and triggers a detector response. The non-robust metric will continue to trigger angular refinement along this one path, but this will not significantly decrease the error in the detector as the four other streaming paths contribute. This highlights the difficulty in identifying when a non-robust goal-based metric is missing streaming paths; clearly we cannot rely on a constant NDOFs with increasing adapt steps (like we could in Figures \ref{sec:2D void problem} \& \ref{sec:2D void problem - 100}, where the adaptivity does not trigger and the NDOFs stays constant). 
\begin{figure}[th]
\centering
\subfloat[][Error vs CDOFs]{\label{fig:3D_scat_box_convg}\includegraphics[width =0.47\textwidth]{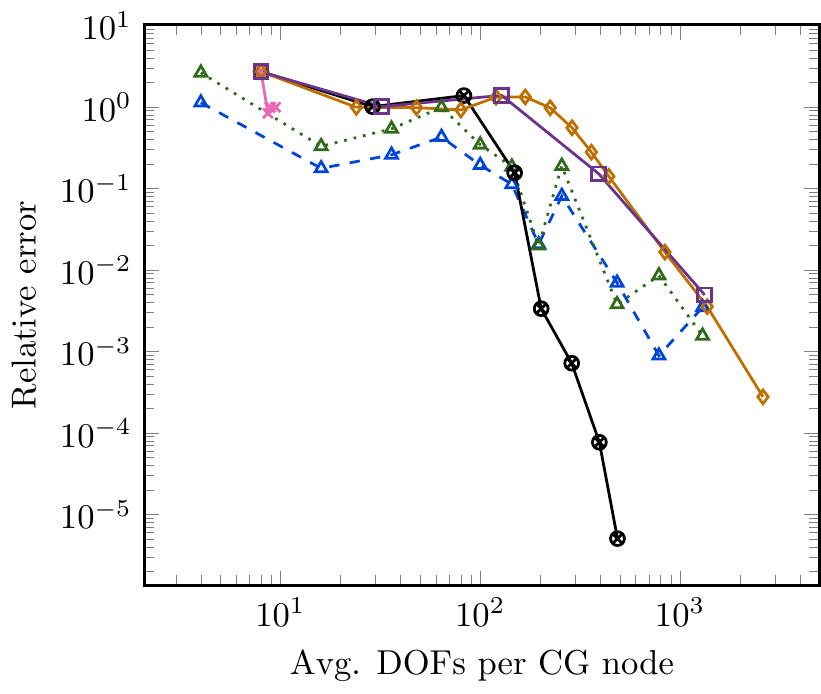}}
\subfloat[][Error vs total runtime]{\label{fig:3D_scat_box_time}\includegraphics[width =0.47\textwidth]{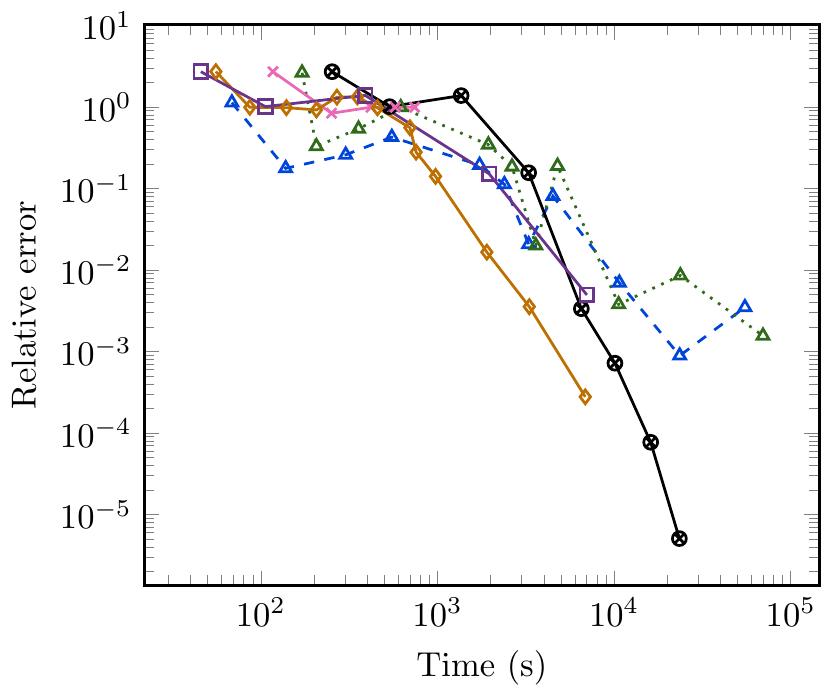}}
\caption{Comparison of the relative error of the detector response, for different angular discretisations for the 3D scatter box problem. The dashed \textcolor{matlabblue}{$\triangle$} is uniform FP$_n$ with $\Sigma_{\textrm{f}}=1$, the dotted \textcolor{foliagegreen}{$\triangle$} is uniform FP$_n$ with $\Sigma_{\textrm{f}}=0.1$, \textcolor{deludedorange}{$\diamond$} uniform LS P$^0$ FEM, the \textcolor{gaylordpurple}{\Square} is non-standard Haar wavelets with fixed angular refinement between $\mu \in [-1, 1]$ and $\omega \in [0, 3.15]$, the \textcolor{darksalmon}{x} are goal-based adapted non-standard Haar wavelets with error target 1\xten{-2} and the \textcolor{black}{$\otimes$} are goal-based adapted non-standard Haar wavelets with robust error target 1\xten{-3} and with reduced tolerance solves.}
\label{fig:3D_scat_box_result}
\end{figure}

\fref{fig:3D_scat_box_convg} shows that the non-robust Haar adapt after 2 steps does not add DOFs in the one streaming path it has resolved and \fref{fig:3D_scat_box_time} shows the subsequent plateu in error vs runtime. We could decrease the adaptive tolerance to force it to further adapt down the single streaming path, but this does not decrease the error given the other more important streaming paths are not resolved. We can also see in \fref{fig:3D_scat_box_result} that the LS P$^0$ and fixed refinement Haars are converging similarly in this problem and that after approximately 200 DOFs are in the asymptotic regime in this problem. The FP$_n$ solutions are in the asymptotic regime at low order (from FP$_1$) but interestingly are converging non-montonically. 
\begin{figure}[th]
\centering
\subfloat[][Spatial position $x=5.3, y=5.13, z=4.76$. The camera is pointed in the $-y$ direction.]{\label{fig:streaming_path_8}\includegraphics[width =0.47\textwidth]{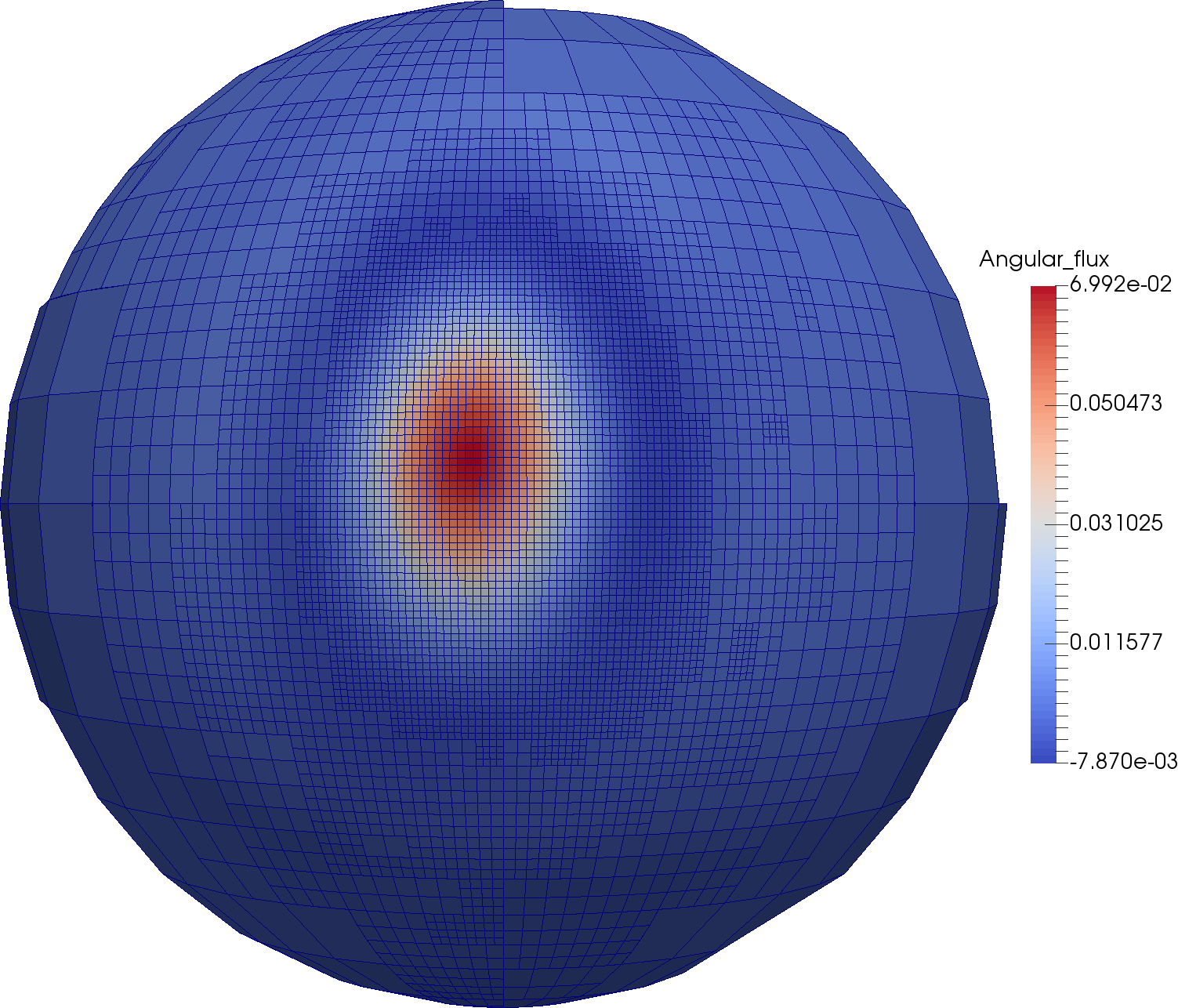}}
\subfloat[][Spatial position $x=0.5, y=5, z=5$. The camera is pointed in the $-x$ direction.]{\label{fig:second_scatter_box}\includegraphics[width =0.47\textwidth]{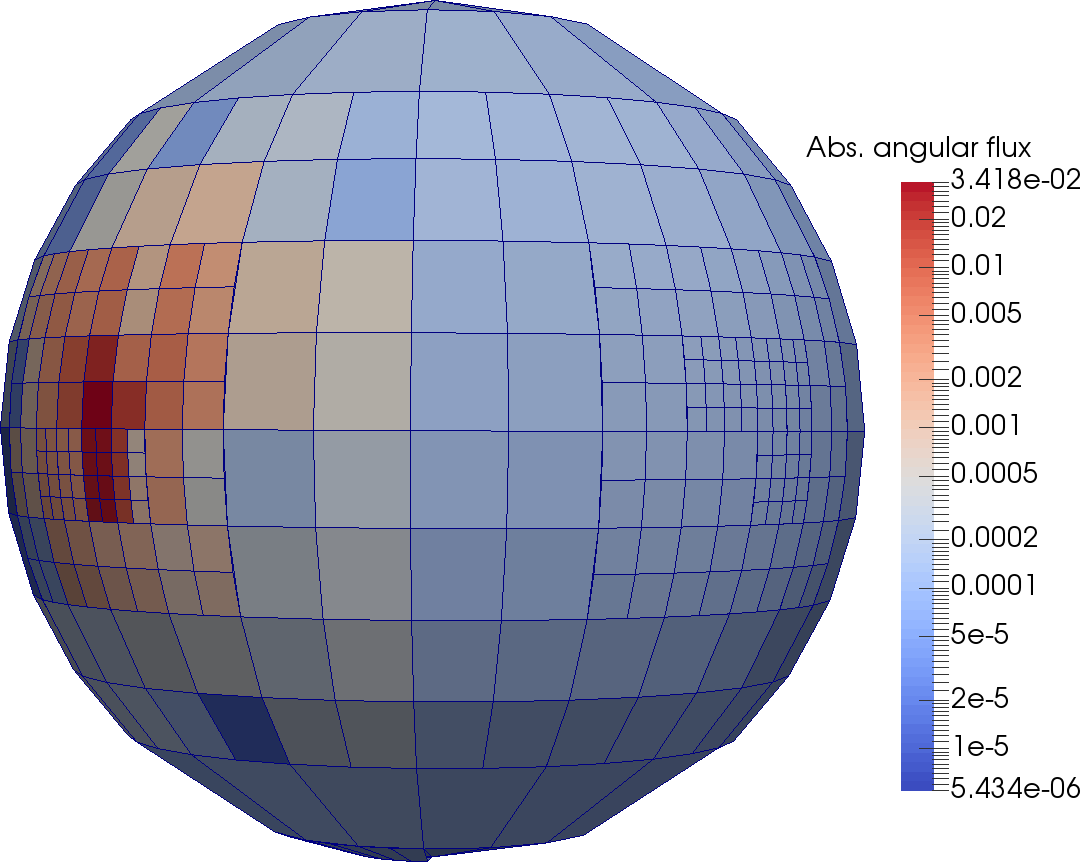}}
\caption{Forward angular flux with adapted wavelets after 8 adapt steps with robust error target 1\xten{-3} at different points in the 3D scatter box problem. The angular discretisations are on the $r=1$ sphere, but have been projected onto faceted polyhedra for ease of visualisation.}
\label{fig:robust_adapt}
\end{figure}

We can see in \fref{fig:3D_scat_box_convg} that the robust adapted Haars are performing very well in this problem, giving a reduction of an order of magnitude in DOFs when compared to the fixed refinement Haars and LS P$^0$ in this problem. We can also see that the FP$_3$ surrogate solution used is sufficient to allow the robust adapt to use less DOFs than the fixed refinement Haars in the pre-asymptotic regime, indicating they are not adapting (even close to) uniformly. The runtime of the robust adapted Haars is also less than the fixed refinement Haar, as shown in \fref{fig:3D_scat_box_time}, and projecting the LS P$^0$ out to further refinement would see the crossover point with the adapt. Computing the forward and adjoint FP$_n$ surrogates only costs 1\% of the total runtime for the highest order refinement with robust Haars. 

The key feature here is that to see a constant reduction in error with each refinement step, the robust adapted Haar must be resolving all the different streaming paths in this problem, not just the one seen by the coarse angular discretisation (and hence the non-robust adapt). This helps verify that our choice of heuristic ratio which determines underresolved nodes is sufficient for this problem, that features different streaming paths of different size/importance. 

\fref{fig:robust_adapt} shows the angular flux in the forward robust adapted Haar solution after 8 adapt steps. \fref{fig:streaming_path_8} in particular shows the flux in the main streaming path between source and detector and we can see heavy anisotropic angular refinement has occured to resolve this important path. \fref{fig:second_scatter_box} however shows the angular flux in the middle of the leftmost scatter box. We can see that refinement only up to level 6 has occured (when a max of 8 levels is possible), and the two angular regions at level 6 are pointing towards the source and detector along the angular equator. This is because the robust metric has identified the indirect streaming path caused by the scattering box, but this streaming path is far less important than the direct path and hence at the adapted tolerance used (1\xten{-2}) it is not the focus of heavy refinement; 6 levels of refinement is enough to resolve this path sufficiently for the given tolerance. Indeed \fref{fig:second_scatter_box} shows that refinement has not occured as might be expected to resolve the path from source, to rightmost scatter box and then finally to the leftmost box; we would expect to see refinement in the bottom two octants visible to the camera. This is because this path is even less important for the tolerance used, given it must scatter twice to reach the detector. Reducing the tolerance would cause both of these paths to be further refined like the direct path (though the direct path will always refine more than the indirects given its importance). 
\begin{table}[ht]
\centering
\begin{tabular}{ l c c c c c c c c c}
\toprule
\textbf{Adapt step:} & \textbf{1} & \textbf{2} & \textbf{3} & \textbf{4} & \textbf{5} & \textbf{6} & \textbf{7} & \textbf{8}\\
\midrule  
Non-robust metric & 4.63 & 1.51 & 0.127 & 0.127 & 0.127 & - & - & -\\
Robust metric & 12.8 & 84.3 & 41.6 & 138.5 & 6280.2 & 32188 & 327326 & 5042647\\
\bottomrule  
\end{tabular}
\caption{Effectivity index for the goal-based adapted discretisation shown in \fref{fig:3D_scat_box_result}, for the 3D scatter box problem.}
\label{tab:3D_scat_box_effec}
\end{table}

Finally \tref{tab:3D_scat_box_effec} shows the effectivity index for both the non-robust and robust Haar adapts. In contrast to the previous two problems, we can see the non-robust metric has a non-zero effectivity index, given the coarse angular discretisation can see one streaming path as discussed. Interestingly, the effectivity index of the robust metric looks worse in the pre-asymptoic region when compared to the non-robust metric. \fref{fig:3D_scat_box_result} of course clearly shows the non-robust metric failing to adapt correctly. This further highlights the difficulty of identifying when the non-robust metric is not resolving streaming paths and helps justify our choice of using the surrogate solutions to drive our adapt, rather than a combination of the non-robust and robust error metrics, as discussed in \secref{sec:Adaptivity algorithm}. \tref{tab:3D_scat_box_effec} also shows the effectivity index growing pathologically in the asymptoic region, as discussed by \cite{Dargaville2019}.
\section{Conclusions}
This paper has presented a method to enable the use of goal-based error metric with angular adaptivity in problems where using a non-rotationally invariant angular discretisations causes severe ray-effects. This method involves computing low-order filtered spherical harmonics forward and adjoint solutions alongside the coarse NRI angular discretisations. This exploits the fact that the pre-asymptotic regime for both angular discretisations is very different. The FP$_n$ solutions are then used to determine which spatial nodes in the problem are considered ``underresolved'' and the FP$_n$ forward and adjoint solutions are then projected into the adapted NRI angular space on those nodes and used within the error metric to allow anisotropic angular adaptivity. Once these spatial nodes are considered resolved the traditional NRI error can take over and adapt to high levels of refinement. 

We tested this method on three simple problems with pure streaming; traditional error metrics would not trigger adaptivity in two of these problems and we found our improved metric matched the error produced by \textit{a priori} refined discretisations. We then deliberately designed a third example problem to trigger incorrect refinement with traditional error metrics and showed our improved metric correctly identified the most important streaming paths in the problem. In all these problems, computing the FP$_n$ forward and adjoint solutions cost at most 1\% of the total runtime of the adaptive process at high order. Furthermore, we verified that the low-order FP$_n$ solutions we used are sufficient to drive adaptivity in these problems.

This is a key point, as computing our FP$_n$ solutions is $\mathcal{O}(n^2)$ as the angular order is increased and hence we are restricted to using low-order if the cost of computing the FP$_n$ solutions is not to outweigh the cost of computing the NRI angular discretisation we are trying to adapt. We found that, as might be expected, if we use too low FP$_n$ order that this causes our wavelet discretisation to ``adapt'' uniformly until it reaches the asymptoic regime. Using higher order allows the pre-asymptotic adapt to be anisotropic, though we found oscillations in the FP$_n$ solution also cause excessive resolution to be applied, though this less resolution than would be applied by uniform refinement. 

Of course we could always construct harder problems that would require very high FP$_n$ order to resolve sufficiently to drive our adapt; our previous work on FP$_n$ adaptivity helps increase the range of problems we could resolve and our future work will investigate allowing both the FP$_n$ and wavelet discretisations to perform simultaneous goal-based angular adapts. For problems beyond this capability, we would argue that this work at least provides a level of robustness in goal-based metrics that was not previously available; traditional goal-based error metrics silently fail in these problems, whereas as long as the FP$_n$ order is increased, the adaptive tolerance is decreased and the ``underresolved'' ratio is decreased, the error metric in this work will trigger refinement correctly. 

We also found that the effectivity index of our error metric was significantly improved in the pre-asymptoic regime by using our FP$_n$ surrogate solutions. The index in the asymptotic regime however must be improved; much like previous work in advection-diffusion problems, we cannot ignore the possible lack of continuity in the true solution and hence the theory our goal-based metrics are based on. This work is designed to tackle the unique problem of ray-effects faced by pre-asymptotic discretisations of the BTE and future work will examine the asymptotic region in more detail.

If we are to have robust error metrics for Boltzmann transport problems, the key message of this paper is that we must test across a range of parameter regimes that include pure streaming, while also ensuring that numerical diffusion, deliberate alignment of geometry and other mitigating factors do not hide the possible failures of our metrics. We believe this work helps extend the range of applicability of goal-based metrics in Boltzmann problems and helps move towards the aim of robust space/angle adaptivity.
\section*{Acknowledgements}
The authors would like to acknowledge the support of the EPSRC through the funding of the EPSRC grant EP/P013198/1.




\section*{References}
\bibliographystyle{model1-num-names}
\bibliography{bib_library}







\end{document}